\pdfoutput=1
\RequirePackage{ifpdf}
\ifpdf 
\documentclass[pdftex]{sigma}
\else
\documentclass{sigma}
\fi

\numberwithin{equation}{section}

\newtheorem{Theorem}{Theorem}[section]
\newtheorem*{Theorem*}{Theorem}
\newtheorem{Corollary}[Theorem]{Corollary}
\newtheorem{Lemma}[Theorem]{Lemma}
\newtheorem{Proposition}[Theorem]{Proposition}
 { \theoremstyle{definition}
\newtheorem{Definition}[Theorem]{Definition}

 }

\begin{document}
\allowdisplaybreaks

\newcommand{\arXivNumber}{2406.16208}

\renewcommand{\PaperNumber}{062}

\FirstPageHeading

\ShortArticleName{Deformation Families of Quasi-Projective Varieties and Symmetric Projective K3 Surfaces}

\ArticleName{Deformation Families of Quasi-Projective Varieties\\ and Symmetric Projective K3 Surfaces}

\Author{Fan XU}

\AuthorNameForHeading{F.~Xu}

\Address{Graduate School of Mathematics, Nagoya University, \\ Furocho, Chikusaku, Nagoya, 464-8602, Japan}
\Email{\href{mailto:fan.xu.e6@math.nagoya-u.ac.jp}{fan.xu.e6@math.nagoya-u.ac.jp}}

\ArticleDates{Received October 31, 2024, in final form July 16, 2025; Published online July 28, 2025}

\Abstract{The main aim of this paper is to construct a complex analytic family of symmetric projective K3 surfaces through a compactifiable deformation family of complete quasi-projective varieties from \smash{$\operatorname{CP}^2 \#9\overline{\operatorname{CP}}^2$}. Firstly, for an elliptic curve $C_0$ embedded in $\operatorname{CP}^2$, let~\smash{$S \cong \operatorname{CP}^2 \#9\overline{\operatorname{CP}}^2$} be the blow up of $\operatorname{CP}^2$ at nine points on the image of $C_0$ and $C$ be the strict transform of the image. Then if the normal bundle satisfies the Diophantine condition, a tubular neighborhood of the elliptic curve $C$ can be identified through a toroidal group. Fixing the Diophantine condition, a smooth compactifiable deformation of $S\backslash C$ over a 9-dimensional complex manifold is constructed. Moreover, with an ample line bundle fixed on $S$, complete K\"ahler metrics can be constructed on the quasi-projective variety~$S\backslash C$. So complete K\"ahler metrics are constructed on each quasi-projective variety fiber of the smooth compactifiable deformation family. Then a complex analytic family of symmetric projective K3 surfaces over a 10-dimensional complex manifold is constructed through the smooth compactifiable deformation family of complete quasi-projective varieties and an analogous deformation family.}


\Classification{14J28; 32G05}

\section{Introduction}

This paper is motivated by the gluing construction of K3 surfaces through two open complements of the closures of tubular neighborhoods of elliptic curves embedded in \smash{$\operatorname{CP}^2 \#9\overline{\operatorname{CP}}^2$} in \cite{Takayuki-Takato1}. Moreover, a gluing construction of projective K3 surfaces was presented in \cite{Takayuki-Takato2}.

For a rational elliptic surface $ X $ defined as the blow-up of a projective plane at nine base points for a pencil of cubics, Hans-Joachim Hein presented some complete Calabi--Yau metrics on $ X \backslash D $ with $D$ being a fiber on $X$ in \cite{Hans}. However, the blow up of a projective plane at arbitrary nine points may not be a rational elliptic surface.

Let $ \langle 1,\tau \rangle$ denote a lattice for $\tau\in Y:=\{\tau \mid \tau\in \mathbb{C},\, \operatorname{\operatorname{Im}}\tau>0 \}$. Then $C_0(\tau)=\mathbb{C}/ \langle 1,\tau \rangle$ is an elliptic curve for $\tau\in Y$. Let $C_0=\mathbb{C}/ \langle 1,\tau_0 \rangle$ with $\tau_0 \in Y$. The smooth elliptic curve $C_0 (\tau)$ can be holomorphically embedded in $\operatorname{CP}^2$ for $\tau \in Y$. The embedding map is induced by the Weierstrass $\wp$-function
\[\wp(z) =\sum _{\lambda \in \langle 1,\tau \rangle \backslash \{0\}} \left(\frac{1}{(z-\lambda )^2}-\frac{1}{\lambda ^2}\right)+\frac{1}{z^2}\]
and the Eisenstein series
\[G_{2k} ( \langle 1,\tau \rangle)=\sum _{\lambda \in \langle 1,\tau \rangle \backslash \{0\}} \lambda ^{-2k}.\]
Here, for $\tau\in Y$, the holomorphic embedding map $f_\tau\colon C_0 (\tau) \rightarrow\operatorname{CP}^2$ is defined by
\[f_\tau([z])=[\wp(z):\wp'(z):1]=[z_1:z_2:z_3]\]
for $[z]\in C_0 (\tau)$ with $[z]\neq[0]$ and $f_\tau([0])=[0:1:0]$. Furthermore,
\[f_\tau(C_0(\tau))=\bigl\{ [z_1:z_2:z_3]|-4 z_1^3+60 G_4 ( \langle 1,\tau \rangle)z_3^2 z_1+140 G_6 ( \langle 1,\tau \rangle)z_3^3+z_2^2 z_3=0\bigr\}\subset\operatorname{CP}^2\]
for $\tau \in Y$.

So $f_\tau(C_0(\tau))$ is a submanifold of $S_0=\operatorname{CP}^2$. Let \smash{$S(\tau)\cong\operatorname{CP}^2 \#9\overline{\operatorname{CP}}^2$} be the blow up of $S_0$ at nine points in the set $Z:=\{p_1,p_2,\dots ,p_9\}$ on $f_\tau(C_0(\tau))$. Now let $C(\tau)$ denote the strict transform of $f_\tau(C_0(\tau))$ in $S(\tau)$. Then $S(\tau)\setminus C(\tau)$ is a quasi-projective variety for $\tau\in Y$. In particular, let $S=S(\tau_0)$ and $C=C(\tau_0)$. Here, the normal bundle of $C(\tau)$ in $S(\tau)$ is assumed to satisfy the Diophantine condition corresponding to a Diophantine number pair proposed in~\cite{Takayuki-Takato2}. Then a complex analytic family of \smash{$\operatorname{CP}^2 \#9\overline{\operatorname{CP}}^2$} over a 9-dimensional complex manifold is constructed exactly. The simple description of this deformation family can be found in \cite{Takayuki-Takato1}.

Elizabeth Gasparim and Francisco Rubilar presented new definitions for deformation family in \cite{Elizabeth-Francisco}, especially for the deformation family of open manifolds. Taking use of the new definitions, the main theorem is as follows.

\begin{Theorem} \label{Theorem1.1}
There is a smooth compactifiable deformation of $S\backslash C$ over a $9$-dimensional complex manifold $T$. In addition, the deformation is differentially trivial along $T$.
\end{Theorem}

As it was proved in \cite{Anna}, $\forall \tau_1, \tau_2 \in\! Y$, if there was an algebraic isomorphism between $S(\tau_1)\backslash C(\tau_1)$ and $S(\tau_2)\backslash C(\tau_2)$, then there was an induced birational morphism between $S(\tau_1)$ and $S(\tau_2)$ which made the square-zero elliptic curves $C(\tau_1)$ and $C(\tau_2)$ isomorphic to each other. Therefore, ${\forall \tau_1 \neq \tau_2 \in Y}$, if $C(\tau_1)$ is not isomorphic to $C(\tau_2)$, then $S(\tau_1)\backslash C(\tau_1)$ is not algebraically isomorphic to $S(\tau_2)\backslash C(\tau_2)$. So through suitable choice of $T$, every two different fibers of the compactifiable deformation of $S\backslash C$ in Theorem~\ref{Theorem1.1} are not algebraically isomorphic to each other.

A K3 surface with an involution was constructed as an example in \cite[Section 7.1.2]{Takayuki-Takato1} by Takayuki Koike and Takato Uehara. Then through fixing suitable ample line bundles on $S(\tau)$, a gluing construction of symmetric projective K3 surfaces will be introduced here.

Firstly, there are ample line bundles properly defined on each fiber of the deformation family of \smash{$\operatorname{CP}^2 \#9\overline{\operatorname{CP}}^2$} \cite{Tomasz-Halszka}.

Secondly, in \cite{Takayuki-Takato2}, there was a simple description for the construction of some complete K\"ahler metrics on $\text{S}(\tau)\backslash \text{C}(\tau)$ corresponding to an ample line bundle on $S(\tau)$. In this paper, through modifying the method proposed in \cite{Takayuki-Takato2}, the analogous complete K\"ahler metrics corresponding to the ample line bundles selected are constructed on each fiber of the smooth compactifiable deformation family of $S \backslash C$.

Then the smooth compactifiable deformation family is proved to be a deformation family of complete quasi-projective varieties.

In particular, Yoshio Fujimoto proposed a condition for the nine blowing up points such that~$S(\tau)$ was a rational elliptic surface in \cite{Yoshio}. Taking use of the same construction method as Section~\ref{sec3} of this paper, we may construct another compactifiable deformation family and each fiber is an open Calabi--Yau manifold. The Calabi--Yau metrics on the fibers were constructed in~\cite{Hans} by Hans-Joachim Hein. However, this situation will not be further discussed in this paper.

Thirdly, a symmetric projective K3 surface can be constructed through gluing the open completement of closure of the tubular neighborhood of $C(\tau)$ in $S(\tau)$ and an analogous manifold.

Then a deformation family of symmetric projective K3 surfaces is constructed as follows.

\begin{Corollary} \label{Corollary1.2}
A deformation family of symmetric projective K3 surfaces over a $10$-dimensional complex manifold along with a symmetric K\"ahler metric on each fiber can be constructed through a smooth compactifiable deformation of $S \backslash C$ and an analogous deformation family.
\end{Corollary}

Moreover, the complex analytic family of symmetric projective K3 surfaces has an injective Kodaira--Spencer map by \cite[Theorem~1.1]{Takayuki-Takato1}.

The following is to introduce the main contents in the following sections of this paper.

In Section~\ref{sec2}, the main goal is to present a gluing construction of symmetric projective K3 surfaces and the toroidal groups corresponding to Diophantine number pairs.

Firstly, the Diophantine condition is introduced for the construction of a tubular neighborhood of an elliptic curve embedded in \smash{$\operatorname{CP}^2 \#9 \overline{\operatorname{CP}}^2$}.

Secondly, fixing a suitable ample line bundle on $S(\tau)$, the open completement of the closure of the tubular neighborhood of $C(\tau)$ in $S(\tau)$ and an analogous open manifold can be glued to a~symmetric projective K3 surface through a suitable map.

Finally, the specific defined toroidal groups related to Diophantine number pairs are just properly defined in a special class of toroidal groups described in \cite[Theorem]{Christian}.

In Section~\ref{sec3}, the main goal is to prove Theorem~\ref{Theorem1.1} and construct complete K\"ahler metrics on each quasi-projective variety fiber of the deformation family through fixing an ample line bundle on $S(\tau)$.

Moreover, through the process of constructing the complete K\"ahler metrics on $S(\tau)\backslash C(\tau)$, the condition for the gluing construction of symmetric projective K3 surfaces can be confirmed.

In Section~\ref{sec4}, the main goal is to prove Corollary~\ref{Corollary1.2}.

\section{Symmetric projective K3 surfaces and toroidal groups}\label{sec2}

\subsection{A gluing construction of symmetric projective K3 surfaces}

Firstly, it is to give the details for a gluing construction of symmetric projective K3 surfaces.

Julius Ross and David Witt Nystr\"om stated that there was a canonical smooth tubular neighborhood for a compact complex submanifold of a K\"ahler manifold which in general would not be holomorphic in \cite{Julius-David}.

V.I.~Arnold proposed a theorem that in most situations an elliptic curve holomorphically embedded in a complex surface with self-intersection index 0 should have a small neighborhood biholomorphically equivalent to a neighborhood of zero section of the normal bundle in \cite{Vladimir}. Ueda also mentioned a corresponding theorem in \cite{Tetsuo}.

Let $\mathbb{N} \subset \mathbb{Z} $ denote the set consisting of positive integers.

Another fact is that any topologically trivial holomorphic line bundle on an elliptic curve is flat \cite{Tetsuo}. The following is a smooth version definition for the monodromy of a complex line bundle over a Riemann surface.

\begin{Definition} \label{Definition2.1}
Let $\nabla$ be a flat connection on an arbitrary complex line bundle $L$ representing an isomorphism class of complex line bundles over a Riemann surface $M$. Define $g_0={\rm e}^{\sqrt{-1} t}\in S^1$ for $t\in [0,2 \pi ]$ and let $\gamma \colon S^1\to M$ be a smooth closed curve. In addition, let $\psi$ denote a parallel section of $\widehat{\gamma}^*L$ with $\widehat{\gamma}=\gamma \circ g_0 $. Then $\widehat{f}(\gamma)\in S^1$ with $\psi_{2\pi} =\widehat{f}(\gamma)\psi_0$ is defined to be the monodromy of $L$ along $\gamma$.
\end{Definition}

The Diophantine condition for a pair of real numbers from \cite{Takayuki-Takato2} is defined as follows.

\begin{Definition}[{\cite[Definition~2.1]{Takayuki-Takato2}}] \label{Definition2.2}
A pair of numbers $(p,q)\in \mathbb{R}^2$ is said to satisfy the Diophantine condition if there exist $\vartheta >0$ and $A>0$ such that
\[
\underset{\mu_0 ,\nu_0 \in \mathbb{Z}}{\min } \big| n\bigl(p+q\sqrt{-1} \bigr)-\bigl(\mu_0 +\nu_0\sqrt{-1} \bigr)\big| \geq A\cdot n^{-\vartheta }
\]
for any $n \in \mathbb{N}$.
\end{Definition}

With the definitions above, the Diophantine condition for a flat holomorphic line bundle on an elliptic curve embedded in a complex surface for the construction of a holomorphic tubular neighborhood of the elliptic curve is defined as follows.

\begin{Definition}[{\cite[Definition~2.2]{Takayuki-Takato2}}] \label{Definition2.3}
For $\tau\in Y$, let $\widetilde{\alpha}$ and $\widetilde{\beta}$ be smooth loops on the elliptic curve $C_0(\tau)=\mathbb{C}/ \langle 1,\tau \rangle$ corresponding to the line segments $[0,1]$ and $[0,\tau]$, respectively. For any topologically trivial line bundle $L\in \operatorname{Pic}^0(C_0(\tau))$, let \smash{${\rm e}^{\widetilde{p}\cdot2\pi\sqrt{-1}}$} and \smash{${\rm e}^{\widetilde{q}\cdot2\pi\sqrt{-1}}$} denote the monodromies of $L$ along the loops $\widetilde{\alpha}$ and \smash{$\widetilde{\beta}$} respectively with $(\widetilde{p},\widetilde{q})\in \mathbb{R}^2$. Then $L$ is said to satisfy the Diophantion condition if $(\widetilde{p},\widetilde{q})$ satisfies the Diophantine condition.
\end{Definition}

According to \cite[Theorem~1.6]{Takayuki-Takato1}, if $\widetilde{L}=N_{C(\tau)/S(\tau)}\in \operatorname{Pic}^0(C(\tau))$ satisfies the Diophantine condition in Definition~\ref{Definition2.3}, \smash{$-\log d\bigl(\mathbb{I}_{C(\tau)},\widetilde{L}^n\bigr)=O(\log n)$} as $n\rightarrow\infty$ with $d$ being an invariant distance on $\operatorname{Pic}^0(C(\tau))$ and $\mathbb{I}_{C(\tau)}$ being the holomorphically trivial line bundle on $C(\tau)$. So from~\cite[Theorem 3]{Tetsuo} which is a generalized result originating from \cite{Vladimir} by V.I.~Arnold, the following theorem is true.

\begin{Theorem}[\cite{Takayuki-Takato2}] \label{Theorem2.4}
If the normal bundle \emph{$N_{C(\tau)/S(\tau)}\in \operatorname{Pic}^0(C(\tau))$} satisfies the Diophantine condition, then $C(\tau)$ admits a holomorphic tubular neighborhood in \smash{\emph{$S(\tau)\cong\operatorname{CP}^2\text{$\#$9}\overline{\operatorname{CP}}^2$}} biholomorphic to a neighborhood of the zero section in \smash{$N_{C(\tau)/S(\tau)}$}.
\end{Theorem}

Now assume that \smash{$N_{C(\tau)/S(\tau)}$} satisfies the Diophantine condition. Let $\alpha$ and $\beta$ be smooth loops on $C(\tau)$ corresponding to the line segments $[0,1]$ and $[0,\tau]$ on $C_0(\tau)$. Let \smash{${\rm e}^{p\cdot2\pi\sqrt{-1}}$} and~\smash{${\rm e}^{q\cdot2\pi\sqrt{-1}}$} be the monodromies of $N_{C(\tau)/S(\tau)}$ along the loops $\alpha$ and $\beta$ respectively with $(p,q)\in \mathbb{R}^2$. So $(p,q)$ should satisfy the Diophantine condition. Taking use of Theorem~\ref{Theorem2.4}, the following statement is true.

There is a tubular neighborhood $W$ of $C(\tau)$ biholomorphic to a neighborhood of the zero section in $N_{C(\tau)/S(\tau)}$ which can be expressed as
$
W \cong \bigl\{(z,w)\in \mathbb{C}^2\mid | w| <r \bigr\}/{\sim}$
for a real number~${r>1}$, where $\sim$ is the equivalence relation generated by
\[(z,w)\sim\bigl(z+1,\exp\bigl(p\cdot 2\pi \sqrt{-1}\bigr)\cdot w\bigr)\sim\bigl(z+\tau,\exp\bigl(q\cdot 2\pi \sqrt{-1}\bigr)\cdot w\bigr)\]
for $\tau\in Y$.

Then for $s\in \Delta\backslash \{0\}$ with $\Delta:=\{s\in \mathbb{C}\mid |s|<1\}$, let
\[M_s:=S(\tau)\backslash\bigl\{[(z,w)]\in W\mid |w| \leq \sqrt{|s|}/r\bigr\}\]
and
\[V_s:=\bigl\{[(z,w)]\in W\mid \sqrt{|s|}/r<|w|<\sqrt{|s|}r\bigr\}.\]

The following is to describe some ample line bundles on $S(\tau)$.

Let $E_i$ be the exceptional divisor corresponding to the point $p_i \in Z:= \{p_1, p_2, \dots , p_9 \}$ for~${i \in \{1, 2, \dots , 9 \}}$, respectively. Moreover, let $\pi\colon S(\tau)\rightarrow \operatorname{CP}^2$ be the blow up of $\operatorname{CP}^2$ at nine points in $Z:= \{p_1, p_2, \dots , p_9 \}\subset f_{\tau}(C_0(\tau))$ and $H=(\pi)^* \mathcal{O}_{\operatorname{CP}^2}(1)$ be the pull back bundle of the hyperplane line bundle. Then the Nakai--Moishezon criterion \cite[Appendix~A]{Robin} can be used to prove the following simplified theorem.

\begin{Lemma}[{\cite[Theorem~3]{Tomasz-Halszka}}] \label{Lemma2.5}
Let $\pi\colon S(\tau)\rightarrow$ \emph{$\operatorname{CP}^2$} be the blowing up of \emph{$\operatorname{CP}^2$} at nine points defined above. For the line bundles of the form $L=d\cdot H-k\cdot\sum _{i=1}^9 E_i$ with $k\geq2$ and $d\geq3k+1$ being integers, if \smash{$9\leq \frac{d^2}{k^2}-1$}, then the line bundles are ample.
\end{Lemma}

Now let $C^{\pm}_0(\tau) \subset S^{\pm}_0=\operatorname{CP}^2 $ be two copies of $f_{\tau}(C_0(\tau))$. Then the map
\[\ell_{\operatorname{CP}^2}([a_1:a_2:a_3])=[a_1:-a_2:a_3] \qquad \text{for} \quad [a_1:a_2:a_3]\in \operatorname{CP}^2\]
is a holomorphic involution on $\operatorname{CP}^2$ and induces a holomorphic involution on $f_{\tau}(C_0(\tau))$ as well as \smash{$C^\pm_0(\tau)$} (i.e., \smash{$(\ell_{\operatorname{CP}^2})^{-1}=\ell_{\operatorname{CP}^2}$}). So there is a biholomorphic map \smash{$\ell_{S^{+}_0}$} induced by \smash{$\ell_{\operatorname{CP}^2}$} from~$S^{+}_0$ to $S^{-}_0$. Since $C^\pm_0(\tau)$ are two copies of $f_{\tau}(C_0(\tau))$, there exists an identity map $g_\tau$ from~$C_0^{+}(\tau)$ to~$C_0^{-}(\tau)$. Let \smash{$Z^{\pm}=\bigl\{p^{\pm}_1, p^{\pm}_2, \dots , p^{\pm}_9\bigr\}$} with \smash{$p^{-}_i=\ell_{S^{+}_0}\bigl(p^{+}_i\bigr)$} for $i \in \{1, 2, \dots , 9 \}$ be the sets collecting the nine blow up points on \smash{$C^\pm_0(\tau)$}, respectively. Moreover, let $p^{+}_0$ and~\smash{$p^{-}_0=\ell_{S^{+}_0}\bigl(p^{+}_0\bigr)$} be inflection points of \smash{$C_0^{\pm}(\tau)$}, respectively. Here, the set \smash{$Z^+ \subset C^+_0(\tau)$} is just a~copy of the set~$Z \subset f_{\tau}(C_0(\tau))$.

In addition, let \smash{$S^{\pm}(\tau)\cong\operatorname{CP}^2 \#9\overline{\operatorname{CP}}^2$} be the blow-ups of $S^{\pm}_0$ at nine points in the sets $Z^{\pm}:=\bigl\{p^{\pm}_1, p^{\pm}_2, \dots , p^{\pm}_9\bigr\}$ on $C_0^{\pm}(\tau)$ and $C^{\pm}(\tau)$ be the strict transforms of $C_0^{\pm}(\tau)$, respectively. Furthermore, let \smash{$\ell_{S^{+}(\tau)}\colon S^{+}(\tau)\rightarrow S^{-}(\tau)$} be the biholomorphic map defined on $S^{+}(\tau)$ induced by the biholomorphic map \smash{$\ell_{S^{+}_0}\colon S^{+}_0\rightarrow S^{-}_0$}.

Then the biholomorphic map $\widetilde{g}_{\tau}\colon C^{+}(\tau)\rightarrow C^{-}(\tau)$ induced by $g_{\tau}$ will map the normal bundle~${N_{C^+(\tau)/S^+(\tau)}\in \operatorname{Pic}^{0}\bigl(C^{+}(\tau)\bigr)}$ to the dual of the normal bundle $N_{C^-(\tau)/S^-(\tau)}\in \operatorname{Pic}^{0}(C^{-}(\tau))$ through $(\widetilde{g}_{\tau})^*N_{C^-(\tau)/S^-(\tau)}\cong (N_{C^+(\tau)/S^+(\tau)})^{-1}$. This can be derived through computing the monodromies of \smash{$N_{C^{\pm}(\tau)/S^{\pm}(\tau)}$}.

Taking use of Lemma~\ref{Lemma2.5}, let $E^{\pm}_i$ be the exceptional divisors corresponding to the points $p^{\pm}_i \in Z^{\pm}=\bigl\{p^{\pm}_1, p^{\pm}_2, \dots , p^{\pm}_9\bigr\}$ for $i \in \{1, 2, \dots , 9 \}$, respectively. Moreover, let $H^{\pm}$ be the pull back bundles of the hyperplane line bundles on $S_0^{\pm}$ respectively with $H^+$ being a copy of $H$ and $(\ell_{S^{+}(\tau)})^*H^-=H^+$. So there are ample line bundles $L^{\pm}$ on $S^{\pm}(\tau)$ such that $\bigl(L^{+}.C^{+}(\tau)\bigr)=(L^{-}.C^{-}(\tau))$ with $(\ell_{S^{+}(\tau)})^*L^{-}=L^{+}$ as follows.

\begin{Theorem} \label{Theorem2.6}
There exist ample line bundles $L^{\pm}$ on $S^{\pm}(\tau)$ such that $\bigl(L^{+}.C^{+}(\tau)\bigr)\!=\!(L^{-}.C^{-}(\tau))$ with $(\ell_{S^{+}(\tau)})^*L^{-}=L^{+}$.
\end{Theorem}

\begin{proof}
Since \smash{$K^{-1}_{S^{+}(\tau)}=\big[C^{+}(\tau)\big]=3H^{+}-\sum _{i=1}^9 E_i^{+}$} and \smash{$K^{-1}_{S^{-}(\tau)}=[C^{-}(\tau)]=3H^{-}-\sum _{i=1}^9 E_i^{-}$}, \smash{$\bigl(L^{+}.C^{+}(\tau)\bigr)=3d^{+}-9k^{+}$} and \smash{$(L^{-}.C^{-}(\tau))=3d^{-}-9k^{-}$}. Therefore, for $L^{+}=d^{+}\cdot H^{+}-k^{+} \cdot \smash{\sum _{i=1}^9 E_i^{+}}$ and \smash{$L^{-}=d^{-}\cdot H^{-}-k^{-}\cdot\sum _{i=1}^9 E_i^{-}$} with $k^{\pm}\geq2$, $d^{+} \geq 3k^{+}+1 $ and $d^{-} \geq 3k^{-}+1 $ being integers, $\bigl(L^{+}.C^{+}(\tau)\bigr)=(L^{-}.C^{-}(\tau))$ if and only if $d^{+}-3k^{+}=d^{-}-3k^{-}$.

Then for $d^{+}=d^{-}$ and $k^{+}=k^{-}$ satisfying \smash{$10\leq \frac{(d^{+})^2}{(k^{+})^2}\leq \frac{(d^{-})^2}{(k^{-})^2}$} with $k^{+}=k^{-}\geq 2$ and $d^{+}=d^{-}\geq 3k^{+}+1=3k^{-}+1$ being integers,
\[
L^{+}=d^{+}\cdot H^{+}-k^{+}\cdot\sum _{i=1}^9 E_i^{+}\qquad \text{and}\qquad L^{-}=d^{-}\cdot H^{-}-k^{-}\cdot\sum _{i=1}^9 E_i^{-}\]
 are ample line bundles on $S^{+}(\tau)$ and $S^{-}(\tau)$ respectively such that $\bigl(L^{+}.C^{+}(\tau)\bigr)=(L^{-}.C^{-}(\tau))$. In addition, $(\ell_{S^{+}(\tau)})^*L^{-}=L^{+}$. So Theorem~\ref{Theorem2.6} is proved.
\end{proof}

In fact, for an arbitrary ample line bundle $L^+$ on $S^+(\tau)$ and a line bundle $L^-$ on $S^-(\tau)$ defined through $(\ell_{S^{+}(\tau)})^*L^{-}=L^{+}$, it is not so hard to prove that $\bigl(L^{+}.C^{+}(\tau)\bigr)=(L^{-}.C^{-}(\tau))$. So the construction in this paper can be generalized.

Now let $\alpha_{\pm}$ and $\beta_{\pm}$ be smooth loops on $C^{\pm}(\tau)$ corresponding to the line segments $[0,1]$ and $[0,\tau]$ on $C_0(\tau)$, respectively. So \smash{${\rm e}^{{\pm}p\cdot2\pi\sqrt{-1}}$} and \smash{${\rm e}^{{\pm}q\cdot2\pi\sqrt{-1}}$} are the monodromies of $N_{C^{\pm}(\tau)/S^{\pm}(\tau)}$ along the loops $\alpha_{\pm}$ and $\beta_{\pm}$, respectively with $(p,q) \in \mathbb{R}^2$ being a Diophantine number pair. Taking use of Theorem~\ref{Theorem2.4}, the following statement is true.

There are tubular neighborhoods $W^{\pm}$ of $C^{\pm}(\tau)$ biholomorphic to neighborhoods of the zero section in \smash{$N_{C^{\pm}(\tau)/S^{\pm}(\tau)}$} such that $W^+$ is just a copy of $W$ and \smash{$W^- =\ell_{S^+(\tau)}\bigl(W^+\bigr)$} for $\tau\in Y$.

Then for arbitrary ample line bundles $L^{\pm}\rightarrow S^{\pm}(\tau)$ satisfying $\bigl(L^{+}.C^{+}(\tau)\bigr)=(L^{-}.C^{-}(\tau))$ with $(\ell_{S^{+}(\tau)})^*L^{-}=L^{+}$, let
\[p^{\pm}=3f^{-1}_{\tau}\bigl(p_0^{\pm}\bigr)H^{\pm}-\sum_{i=1}^9f^{-1}_{\tau}\bigl(p_i^{\pm}\bigr)E^{\pm}_i\]
and $b_0=\bigl(L^{+}.C^{+}(\tau)\bigr)=(L^{-}.C^{-}(\tau))$.Then there exists a unique $\xi=\frac{1}{b_0}\bigl((p^-.L^-)-\bigl(p^+.L^+\bigr)\bigr)\in\mathbb{C}$ up to modulo $\langle1,\tau\rangle$ such that $g_{\xi}=\ell_{\xi}\circ\widetilde{g}_{\tau}$ with $\ell_{\xi}$ being a translation induced from $\mathbb{C}\ni z\rightarrow z+\xi\in\mathbb{C}$ gives the following gluing construction \cite[Proposition~5.1]{Takayuki-Takato2}. In our case, $\xi =0$ and therefore $g_{\xi}$ coincides with $\widetilde{g}_{\tau}$.

For $s\in \Delta \backslash \{0\}$ with $\Delta:=\{s\in \mathbb{C}\mid |s|<1\}$, let $M_s^+\subset S^+(\tau)$ be a copy of $M_s$, $M^-_s =\ell_{S^+(\tau)}\bigl(M^+_s\bigr) \subset S^-(\tau)$, $V^+_s \subset S^+(\tau)$ be a copy of $V_s$ and $V^-_s=\ell_{S^+(\tau)}\bigl(V^+_s\bigr) \subset S^-(\tau)$ with
$f_s\colon V^{+}_s \rightarrow V^{-}_s$ defined by
\[f_s\bigl(\big[\bigl(z^+,w^+\bigr)\big]\bigr)=\big[\bigl(g_{\xi}\bigl(z^+\bigr),s/w^+\bigr)\big] \qquad \text{for} \quad \big[\bigl(z^+,w^+\bigr)\big]\in V^+_s.\]

The following is the definition for the symmetric K3 surfaces here.

\begin{Definition} \label{Definition2.7}
For a K3 surface $\mathcal{K}$, if there exists an ample line bundle $L_{\mathcal{K}}$ together with a~non-trivial holomorphic involution $f_{\mathcal{K}}$ defined on $\mathcal{K}$ such that $f_{\mathcal{K}}^*L_{\mathcal{K}}=L_{\mathcal{K}}$, then $\mathcal{K}$ is said to be symmetric with respect to $L_{\mathcal{K}}$.
\end{Definition}

So the gluing construction of symmetric K3 surfaces is as follows.

\begin{Theorem} \label{Theorem2.8}
There exists a~sufficiently small $\varepsilon_0 >0$ such that for any $s \in \Delta \setminus \{ 0\}$ and ${0<|s|<\varepsilon_0}$, identifying $V^{+}_s$ and $V^{-}_s$ through the biholomorphic map $f_s$ defined above, the two open complements $M^{\pm}_s$ of the closures of the tubular neighborhoods $W^{\pm}$ of $C^{\pm}(\tau)$ are glued to be a~symmetric projective K3 surface $X_s$.
\end{Theorem}

\begin{proof}
Let $\ell_{S^+(\tau)\setminus C^+(\tau)}$ be defined as $\ell_{S^+(\tau)\setminus C^+(\tau)}=\ell_{S^+(\tau)}|_{S^+(\tau)\setminus C^+(\tau)}$ for $\tau \in Y$. Since $\ell_{S^+(\tau)}\colon S^+(\tau)\rightarrow S^-(\tau)$ is a biholomorphic map with $\ell_{S^+(\tau)}|_{C^{+}(\tau)}\colon C^+(\tau)\rightarrow C^-(\tau)$ being biholomorphic, taking use of \cite[Proposition~2.1]{Takayuki-Takato1}, there are global holomorphic non-vanishing 2-forms~$\eta^{\pm}_s$ defined on $S^{\pm}(\tau)\setminus C^{\pm}(\tau)$ with \smash{$\ell^{*}_{S^{+}(\tau)\setminus C^{+}(\tau)}(\eta^{-}_s)=\eta^{+}_s$} such that the restrictions of~$\eta^{\pm}_s$ on two open complements $M^{\pm}_s$ of the closures of the tubular neighborhoods of $C^{\pm}(\tau)$ can be glued through $f_s$ to be a global holomorphic non-vanishing 2-form $\sigma_s$ on $X_s$. In addition,
\[\sigma_s|_{V^+_s}=a\cdot \frac{{\rm d}z^+ \wedge {\rm d}w^+}{w^+}\]
for $a \in \mathbb{C} \setminus \{ 0 \}$ . So $\sigma_s$ can be assumed to be normalized in this paper. That is to say, $a=1$. Then $X_s$ is a K3 surface.

\begin{Lemma}[{\cite[Theorem~1.2]{Takayuki-Takato2}}]\label{Lemma2.9} Let $L_s$ be the holomorphic line bundle on $X_s$ derived from $L^{\pm}$ through the gluing construction above. Then $\exists$ sufficiently small $\varepsilon_0 >0$ such that $L_s$ is ample for any $s \in \Delta \setminus \{ 0\} $ and $0<|s|<\varepsilon_0$.
\end{Lemma}

Taking use of Lemma~\ref{Lemma2.9}, for any $s \in \Delta \setminus \{ 0\} $ and $0<|s|<\varepsilon_0$, $L_s$ is an ample line bundle on $X_s$. So $X_s$ is projective for $s \in \Delta \setminus \{ 0\} $ and $0<|s|<\varepsilon_0$.

Let $F_{X_s}\colon X_s \rightarrow X_s$ be defined by
\[F_{X_s}(x)=\begin{cases}
\ell_{M^+_s}(x), &x \in M^{+}_s,\\
\ell^{-1}_{M^+_s}(x), &x \in M^{-}_s
\end{cases}\]
with \smash{$\ell_{M^{+}_s}$} being the restriction of $\ell_{S^{+}(\tau)}$ on $M^{+}_s$. Then $F_{X_s}$ is a non-trivial holomorphic involution on the K3 surface $X_s$. In addition, since $(\ell_{S^{+}(\tau)})^*L^{-}=L^{+}$ , $F^{*}_{X_s}L_s=L_s$ for $s \in \Delta \setminus \{ 0\} $ and~${0<|s|<\varepsilon_0}$. So $X_s$ is a symmetric projective K3 surface for $s \in \Delta \setminus \{ 0\} $ and $0<|s|<\varepsilon_0$. Then Theorem~\ref{Theorem2.8} is proved.
\end{proof}

Here, the existence of $\varepsilon_0$ will also be confirmed in Section~\ref{sec3} during the process of the construction of complete K\"ahler metrics on the quasi-projective varieties $S(\tau)\setminus C(\tau)$ for $\tau \in Y$.

\subsection{Toroidal groups and tubular neighborhoods}

The following is to introduce the corresponding toroidal groups $U_0=\mathbb{C}^2/\Lambda_0$ with
$\Lambda_0=\big\langle\big(\begin{smallmatrix}0 \\1\end{smallmatrix}\big),\allowbreak\big(\begin{smallmatrix}1 \\p\end{smallmatrix}\big),\big(\begin{smallmatrix}\tau \\q\end{smallmatrix}\big)\big\rangle$
provided in \cite{Takayuki-Takato2} for the Diophantine pair $(p,q)$ and show that $U_0=\mathbb{C}^2/\Lambda_0$ is a~quasi-abelian variety.

\begin{Definition}[{\cite[Definition~1.1.1]{Yukitaka-Klaus}}] \label{Definition2.10}
 Let $\Lambda$ be a lattice which is a discrete subgroup of $\mathbb{C}^n$ for $n\in \mathbb{N}$. Then $\mathbb{C}^n/\Lambda$ is called a toroidal group if there does not exist any non-constant holomorphic function on $\mathbb{C}^n/\Lambda$.
\end{Definition}

The toroidal group is a topological group. Moreover, it is an abelian complex Lie group~\cite{Yukitaka-Klaus}.

\begin{Lemma}[{\cite[Theorem~1.1.4]{Yukitaka-Klaus}}] \label{Lemma2.11}
 For $\mathbb{C}^2/\Lambda$ with $\Lambda$ a discrete subgroup of $\mathbb{C}^2$, $\mathbb{C}^2/\Lambda$ is a~toroidal group if and only if there does not exist any~${\sigma \in \mathbb{C}^2 \backslash \{ 0 \}}$ such that the scalar product~${\langle \sigma,\lambda \rangle \in \mathbb{Z}}$ is integral for all $ \lambda \in \Lambda $.
\end{Lemma}

Taking use of Lemma~\ref{Lemma2.11}, the following lemma can be derived. Here, let $\mathbb{Q}$ denote the set consisting of rational numbers.

\begin{Lemma} \label{Lemma2.12}
For any real number pair $(p,q)$ satisfying the Diophantine condition and complex number $\tau \in Y$, $U_0=\mathbb{C}^2/\Lambda_0$ with $\Lambda_0=\bigl\langle\bigl(\begin{smallmatrix}0 \\1\end{smallmatrix}\bigr),\bigl(\begin{smallmatrix}1 \\p\end{smallmatrix}\bigr),\bigl(\begin{smallmatrix}\tau \\q\end{smallmatrix}\bigr)\bigr\rangle$ is a toroidal group.
\end{Lemma}

\begin{proof}
For $\sigma =\bigl(\begin{smallmatrix}\sigma_1\\\sigma_2\end{smallmatrix}\bigr) \in \mathbb{C}^2 \backslash \{ 0 \}$,
\begin{gather*}
\left\langle \sigma,\begin{pmatrix}0 \\1\end{pmatrix} \right\rangle =\left\langle \begin{pmatrix}\sigma_1\\
\sigma_2\end{pmatrix}, \begin{pmatrix}0 \\1\end{pmatrix} \right\rangle =\sigma_2, \qquad \left\langle \sigma,\begin{pmatrix}1 \\p\end{pmatrix} \right\rangle=\left\langle \begin{pmatrix}\sigma_1\\
\sigma_2\end{pmatrix}, \begin{pmatrix} 1 \\p \end{pmatrix} \right\rangle =\sigma_1+p\sigma_2
\end{gather*}
and $\bigl\langle \sigma\bigl(\begin{smallmatrix}\tau \\q\end{smallmatrix}\bigr) \bigr\rangle =\bigl\langle \bigl(\begin{smallmatrix}\sigma_1\\
\sigma_2\end{smallmatrix}\bigr), \bigl(\begin{smallmatrix} \tau \\ q \end{smallmatrix}\bigr) \bigr\rangle =\sigma_1 \tau+q\sigma_2$.

If $q \in \mathbb{Q}$, since $(p,q) \in \mathbb{R}^2$ satisfies the Diophantine condition, $p \notin \mathbb{Q}$. Under the condition that $\sigma_2 \in \mathbb{Z}$ and $\sigma_1 \tau + q\sigma_2\in \mathbb{Z}$, then $\sigma_1=0$ or $\sigma_1 \in \mathbb{C} \setminus \mathbb{R}$. So $\sigma_1 + p\sigma_2\notin \mathbb{Z}$.

If $q \notin \mathbb{Q}$, under the condition that $\sigma_2 \in \mathbb{Z}$ and $\sigma_1 \tau + q\sigma_2\in \mathbb{Z}$, then $\sigma_1 \in \mathbb{C} \setminus \mathbb{R}$. So $\sigma_1 + p\sigma_2\notin \mathbb{Z}$.

Therefore, there does not exist any $\sigma \in \mathbb{C}^2 \backslash \{ 0 \}$ such that the scalar product $\langle \sigma,\lambda \rangle \in \mathbb{Z}$ is integral for all $ \lambda \in \Lambda_0 $.

Taking use of Lemma~\ref{Lemma2.11}, $U_0=\mathbb{C}^2/\Lambda_0$ with $\Lambda_0=\bigl\langle\bigl(\begin{smallmatrix}0 \\1\end{smallmatrix}\bigr),\bigl(\begin{smallmatrix}1 \\p\end{smallmatrix}\bigr),\bigl(\begin{smallmatrix}\tau \\q\end{smallmatrix}\bigr)\bigr\rangle$ is a toroidal group for any $(p,q) \in \mathbb{R}^2$ satisfying the Diophantine condition. So Lemma~\ref{Lemma2.12} is proved.
\end{proof}

\begin{Theorem}[{\cite[Proposition~2.4]{Takayuki-Takato2}}] \label{Theorem2.13}
 For $(p,q) \in \mathbb{R}^2$ satisfying the Diophantine condition and $\tau \in Y$, every topologically trivial line bundle on the toroidal group defined as $U_0=\mathbb{C}^2/\Lambda_0$ with $\Lambda_0=\bigl\langle\bigl(\begin{smallmatrix}0 \\1\end{smallmatrix}\bigr),\bigl(\begin{smallmatrix}1 \\p\end{smallmatrix}\bigr),\bigl(\begin{smallmatrix}\tau \\q\end{smallmatrix}\bigr)\bigr\rangle$ is homogeneous.
\end{Theorem}

\begin{proof}
Taking use of the result from \cite{Christian}, there are nine equivalent conditions for the special class of toroidal groups with homogeneous topologically trivial line bundles \cite[Theorem]{Christian}. The following lemma is to show the rewritten condition 9.

\begin{Lemma}[{\cite{Christian}}] \label{Lemma2.14}
 For $(p,q) \in \mathbb{R}^2$ satisfying the Diophantine condition and $\tau \in Y$, let $D_0=\bigl(\begin{matrix}p& q\end{matrix}\bigr)$. Then every topologically trivial line bundle on the toroidal group defined as $U_0=\mathbb{C}^2/\Lambda_0$ with $\Lambda_0=\bigl\langle\bigl(\begin{smallmatrix}0 \\1\end{smallmatrix}\bigr),\bigl(\begin{smallmatrix}1 \\p\end{smallmatrix}\bigr),\bigl(\begin{smallmatrix}\tau \\q\end{smallmatrix}\bigr)\bigr\rangle$ is homogeneous if and only if there exist $ c>0$ and $a \geq 0$ such that~${\|\sigma^t \cdot D_0+ \delta^t\| \geq c{\rm e}^{-a|\sigma|}}$
for all $\sigma \in \mathbb{Z} \setminus \{ 0 \}$ and $\delta \in \mathbb{Z}^2$.
\end{Lemma}

Since $(p,q)\in \mathbb{R}^2$ satisfying the Diophantine condition and $\tau \in Y$, there exist $\vartheta >0$ and~${A>0}$ such that
\[
\underset{\mu_0 ,\nu_0 \in \mathbb{Z}}{\min } \big| n\bigl(p+q\sqrt{-1} \bigr)-\bigl(\mu_0 +\nu_0\sqrt{-1} \bigr)\big| \geq A\cdot n^{-\vartheta }
\]
for any $n \in \mathbb{N}$.

That is to say,
\[\| \sigma^t \cdot D_0+ \delta^t\|=\| \sigma(p\quad q)+ \delta^t\| \geq A\cdot {\rm e}^{-\vartheta|\sigma|}\]
for all $\sigma \in \mathbb{Z} \backslash \{ 0 \}$ and $\delta \in \mathbb{Z}^2$. Therefore, the nine equivalent conditions for the special class of toroidal groups \cite[Theorem]{Christian} are all satisfied here. That is to say, every topologically trivial line bundle on the toroidal group defined as $U_0=\mathbb{C}^2/\Lambda_0$ with $\Lambda_0=\bigl\langle\bigl(\begin{smallmatrix}0 \\1\end{smallmatrix}\bigr),\bigl(\begin{smallmatrix}1 \\p\end{smallmatrix}\bigr),\bigl(\begin{smallmatrix}\tau \\q\end{smallmatrix}\bigr)\bigr\rangle$ is homogeneous. So Theorem~\ref{Theorem2.13} is proved.
\end{proof}

With the above result, from \cite{Christian}, the complex vector space $H^1(U_0,\mathcal{O})$ is finite-dimensional. The following is to show that $U_0=\mathbb{C}^2/\Lambda_0$ defined above are quasi-abelian varieties.

Let
\[\mathbb{R}_{\Lambda_{0}}:=\left \{ x_1 \cdot \begin{pmatrix}0 \\1\end{pmatrix}+x_2 \cdot \begin{pmatrix}1 \\p\end{pmatrix} +x_3 \cdot \begin{pmatrix}\tau \\q\end{pmatrix}\mid x_1\in \mathbb{R},\,x_2\in \mathbb{R},\,x_3\in \mathbb{R}\right \}\]
be the $\mathbb{R}$-span of $\Lambda_0$. Let
\[M\mathbb{C}_{\Lambda_0}:=\mathbb{R}_{\Lambda_0}\cap\sqrt{-1}\mathbb{R}_{\Lambda_0}=\left \{ \begin{pmatrix}x \\0\end{pmatrix}\mid x\in \mathbb{C}\right \}\]
be the maximal $\mathbb{C}$-linear subspace of $\mathbb{R}_{\Lambda_{0}}$.

\begin{Definition}[{\cite[Definition~3.1.6]{Yukitaka-Klaus}}]\label{Definition2.15}
 An ample Riemann form $H_0$ for the discrete subgroup~${\Lambda \subset \mathbb{C}^n}$ of complex rank n is a Hermitian form $H_0$ on $\mathbb{C}^n$ such that
\begin{itemize}\itemsep=0pt
\item[$(1)$]$\operatorname{ Im}H_0|_{\Lambda\times\Lambda}$ is $\mathbb{Z}$-valued,
\item[$(2)$]$H_0$ is positive definite on the maximal $\mathbb{C}$-linear subspace $M\mathbb{C}_\Lambda$ of $\mathbb{R}_\Lambda$.
\end{itemize}
Here, $\mathbb{R}_\Lambda$ is the $\mathbb{R}$-span of $\Lambda$. $M\mathbb{C}_\Lambda$ is the maximal $\mathbb{C}$-linear subspace of $\mathbb{R}_\Lambda$.
\end{Definition}

\begin{Definition}[{\cite[Definition~3.1.6]{Yukitaka-Klaus}}]\label{Definition2.16}
 A quasi-abelian variety is a toroidal group $\mathbb{C}^n/\Lambda$ with an ample Riemann form for the lattice $\Lambda$.
\end{Definition}

Then the following is true.

\begin{Theorem}\label{Theorem2.17}
For $(p,q) \in \mathbb{R}^2$ satisfying the Diophantine condition and $\tau \in Y$, $U_0={\mathbb{C}^2_{(z,\eta)}}/{\Lambda_0}$ with \smash{$\Lambda_0=\bigl\langle\bigl(\begin{smallmatrix}0 \\1\end{smallmatrix}\bigr),\bigl(\begin{smallmatrix}1 \\p\end{smallmatrix}\bigr),\bigl(\begin{smallmatrix}\tau \\q\end{smallmatrix}\bigr)\bigr\rangle$} is a quasi-abelian variety.
\end{Theorem}

\begin{proof}
Let
\[
G=\left(\begin{matrix}\frac{1}{\operatorname{\operatorname{Im}}\tau} & 0\\0 & 0\end{matrix}\right),
\]
 then~$G$ defines a hermitian form on~$\mathbb{C}^2$. $\operatorname{ Im}G|_{\Lambda_0\times \Lambda_0}$ is $\mathbb{Z}$-valued and~$G$ is positive definite on the maximal $\mathbb{C}$-linear subspace $M\mathbb{C}_{\Lambda_0}$ of $\mathbb{R}_{\Lambda_{0}}$. So
 \[
 G=\left(\begin{matrix}\frac{1}{\operatorname{\operatorname{Im}}\tau} & 0\\0 & 0\end{matrix}\right)
 \]
 is an ample Riemann form for the lattice $\Lambda_0$. Therefore, according to Definition~\ref{Definition2.16}, $U_0=\mathbb{C}^2/\Lambda_0$ is a quasi-abelian variety. Then Theorem~\ref{Theorem2.17} is proved.
\end{proof}

Since $M\mathbb{C}_{\Lambda_0}:=\mathbb{R}_{\Lambda_0}\cap\sqrt{-1}\mathbb{R}_{\Lambda_0}=\big \{ \bigl(\begin{smallmatrix}x \\0\end{smallmatrix}\bigr)\mid x\in \mathbb{C}\big \}$ is the maximal $\mathbb{C}$-linear subspace of~$\mathbb{R}_{\Lambda_0}$ and $M\mathbb{C}_{\Lambda_0}$ has the complex dimension 1, the toroidal group $U_0=\mathbb{C}^2/\Lambda_0$ is of type~1 \cite[Defi\-nition~1.2.1]{Yukitaka-Klaus}.

For \smash{$G=\bigl(\begin{smallmatrix}\frac{1}{\operatorname{\operatorname{Im}}\tau} & 0\\0 & 0\end{smallmatrix}\bigr)$}, $\operatorname{\operatorname{Im}}G$ has rank 2 on $\mathbb{R}_{\Lambda_0}$. Moreover, since the toroidal group $U_0=\mathbb{C}^2/\Lambda_0$ is of type 1, the ample Riemann form \smash{$G=\bigl(\begin{smallmatrix}\frac{1}{\operatorname{\operatorname{Im}}\tau} & 0\\0 & 0\end{smallmatrix}\bigr)$} for $\Lambda_0$ is said to be of kind $\frac{2-2\times1}{2}=0$ (see~\cite[Definition~3.1.12]{Yukitaka-Klaus}).

The following theorem will give a fibration of the toroidal group $U_0=\mathbb{C}^2/\Lambda_0$.

\begin{Theorem}[{\cite[Theorem~3.1.16]{Yukitaka-Klaus}}]\label{Theorem2.18}
 For a toroidal group $\mathbb{C}^n/\Lambda$ of type q with $\Lambda$ being a lattice, $\mathbb{C}^n/\Lambda$ is a quasi-abelian variety with an ample Riemann form of kind k if and only if $\mathbb{C}^n/\Lambda$ has a~maximal closed Stein subgroup $K \cong \mathbb{C}^k\times \mathbb{C}^{*m}$ with $2k+m=n-q$ and $(\mathbb{C}^n/\Lambda)/K$ is an abelian variety of dimension $q+k$.
\end{Theorem}

Taking use of Theorem~\ref{Theorem2.18}, the toroidal group $U_0=\mathbb{C}^2/\Lambda_0$ has a maximal closed Stein subgroup $K \cong \mathbb{C}^*$ and $U_0/K$ is an abelian variety of dimension 1.

Furthermore, from \cite{Takayuki-Takato2}, through the map $f_{U_0}\colon U_0=\mathbb{C}^2_{(z,\eta)}/{\Lambda_0}\rightarrow V_{0,\infty}$ defined by sending
\[
[(z,\eta)]\in U_0=\mathbb{C}^2_{(z,\eta)}/{\Lambda_0} \qquad \text{to} \quad \big[\bigl(z,{\rm e}^{2\pi \sqrt{-1}\eta}\bigr)\big]\in V_{0,\infty},
\]
$U_0=\mathbb{C}^2_{(z,\eta)}/{\Lambda_0}\cong V_{0,\infty}$ with $V_{0,\infty}$ defined by
$V_{b_1,b_2}:=\big \{ (z,w)\in \mathbb{C}^2\mid b_1<|w|<b_2 \big \}/{\sim}$
for~${0\leq b_1<b_2\leq\infty}$, where $\sim$ is the equivalence relation generated by
\[(z,w)\sim\big(z+1,\exp\big(p\cdot 2\pi \sqrt{-1}\big)\cdot w)\sim\big(z+\tau,\exp\big(q\cdot 2\pi \sqrt{-1}\big)\cdot w\big)\]
with $\tau \in Y$.

From the definitions above, $V_{0,r}\cup C(\tau)=W$. Let $\pi_W\colon W \rightarrow C(\tau)$ be the natural projection defined as $\pi_W([(z,w)])=[z]$ for $[(z,w)]\in W$. The next step is to give a theorem for the line bundles on the tubular neighborhoods defined above.

\begin{Theorem}[{\cite[Proposition~3.3]{Takayuki-Takato2}}]\label{Theorem2.19}
 For any $L \in\operatorname{ Pic}(W)$, $L$ is the pull back line bundle~$(\pi_W)^*(L|_{C(\tau)})$ of the restricted line bundle $L|_{C(\tau)}$ on $C(\tau)$. Here, the assumption is as above.
\end{Theorem}

This theorem is proved in \cite{Takayuki-Takato2} through the toroidal group defined above.

\section[A deformation family of complete non-compact K\"ahler manifolds]{A deformation family of complete non-compact\\ K\"ahler manifolds}\label{sec3}

\subsection[A deformation family of quasi-projective varieties from CP\string^2 \#9CP\string^2]{A deformation family of quasi-projective varieties from $\boldsymbol{\operatorname{CP}^2 \#9 \overline{\operatorname{CP}}^2}$}\label{sec3.1}

The following is to prove Theorem~\ref{Theorem1.1}. First of all, the definition for a complex analytic family of compact complex manifolds is as follows.

\begin{Definition}[{\cite[Definition~2.8]{Kunihiko}}]\label{Definition3.1}
 For a collection of compact complex manifolds $ {\{M_t \mid t \!\in\! B\}}$ with $B$ being a domain of $\mathbb{C}^m$ and $m\in \mathbb{N} $, $ \{M_t \mid t \in B\}$ is called a complex analytic family of compact complex manifolds if there is a complex manifold $\mathcal{M}$ and a holomorphic map $\pi$ from~$\mathcal{M}$ onto $B$ satisfying the following conditions:
\begin{itemize}\itemsep=0pt
\item[$(1)$] The rank of the Jacobian matrix of $\pi$ is equal to $m$ at every point of $\mathcal{M}$.
\item[$(2)$] For each $t \in B$, $\pi^{-1}(t)=M_t$ is a compact complex submanifold of $\mathcal{M}$.
\end{itemize}
\end{Definition}

Let $(\mathcal{M},B,\pi)$ denote the complex analytic family defined above. Here, a complex analytic family of \smash{$\operatorname{CP}^2 \#9 \overline{\operatorname{CP}}^2$} over a 9-dimensional complex manifold will be presented before constructing a deformation family of quasi-projective varieties.

However, a complex analytic family of elliptic curves over a 1-dimensional complex manifold should be presented at first. Here, for $\tau_0\in Y$, let $U \subset Y$ be a sufficiently small disc neighborhood of $\tau_0$ in $Y$ and $\pi_{\operatorname{CP}^2 \times U}\colon \operatorname{CP}^2\times U\rightarrow U$ be the natural projection map. Let $C_0(\tau)$ be an elliptic curve embedded in
\smash{$\pi_{\operatorname{CP}^2\times U}^{-1}(\tau)=\operatorname{CP}^2 \times\{\tau\}$} as $f_{\tau}(C_0(\tau))\times \{ \tau \}$
for $\tau \in U$ and
\[\widetilde{\mathcal{S}}=\bigl\{ (x, \tau) \in \operatorname{CP}^2 \times U\mid x \in f_{\tau}(C_0(\tau)) \bigr\}.
\]
 Here, $C_0( \tau_0)=C_0$. Then the following lemma is true.

\begin{Lemma}\label{Lemma3.2}
\emph{$\bigl(\widetilde{\mathcal{S}},U,\pi_{\operatorname{CP}^2 \times U}|_{\widetilde{\mathcal{S}}}\bigr)$} is a complex analytic family of elliptic curves over the $1$-di\-men\-sion\-al complex manifold $U$ and \emph{$\bigl(\operatorname{CP}^2\times U, U,\pi_{\operatorname{CP}^2\times U}\bigr)$} is a complex analytic family of \emph{ $\operatorname{CP}^2$} over the $1$-dimensional complex manifold $U$.
\end{Lemma}

The complete proof of Lemma~\ref{Lemma3.2} is in Appendix~\ref{appendixA}.

Now let $U \subset Y$ be a sufficiently small disc neighborhood of $\tau_0$ and $U_\nu$ be sufficiently small disc neighborhoods of $q_\nu$ in the universal cover $\mathbb{C}$ of $ \mathbb{C}/ \langle 1,\tau \rangle=C_0({\tau})$ for $\nu \in \{1,2,3,4,5,6,7,8 \}$. Here, $C_0(\tau_0)=C_0$. Let $f_{\widetilde{\tau},\nu}\colon U_{\nu}\!\rightarrow U_{_\nu}/\langle 1, \widetilde{\tau}\rangle$ be the natural projection for $\nu \in \{1,2,3,4,5,6,7,8 \}$ and $\widetilde{\tau}\in U$.

Let $(p,q)$ still be a Diophantine number pair. For $\hat{t}=(\widetilde{\tau},\hat{p}_1,\dots ,\hat{p}_8)\in T =U \times U_1 \times \dots \times U_8$, let $f_9\bigl(\hat{t}\bigr) \in C_0(\widetilde{\tau})$ be the point fixed by the equation
\[9p_0-\sum _{j=1}^8f_{\widetilde{\tau},j}(\hat{p}_j)-f_9\bigl(\hat{t}\bigr)=q-p\cdot \widetilde{\tau} \quad \mod \langle 1, \widetilde{\tau} \rangle \]
with $p_0$ being an inflection point \cite{Takayuki-Takato2}. In addition, let $\pi\colon\operatorname{CP}^2 \times T \rightarrow T$ be the natural projection and let $\mathcal{S}_0=\widetilde{\mathcal{S}} \times U_1 \times \dots \times U_8$. Then $\bigl(\operatorname{CP}^2 \times T,T, \pi\bigr)$ is a complex analytic family of $\operatorname{CP}^2$.

Through blowing up nine points $ \{ f_{\widetilde{\tau}}(p_1),f_{\widetilde{\tau}}(p_2),\dots ,f_{\widetilde{\tau}}(p_9) \} \times \{ (\widetilde{\tau},\hat{p}_1,\dots ,\hat{p}_8) \}$ on
\[\pi|^{-1}_{\mathcal{S}_0}\bigl(\hat{t}\bigr)=\pi_{\operatorname{CP}^2 \times U}|^{-1}_{\widetilde{S}}(\widetilde{\tau})\times ] \{ (\hat{p}_1,\dots ,\hat{p}_8) ] \}
\]
 with
 \[
 \big \{ p_1=f_{\widetilde{\tau},1}(\hat{p}_1),\dots ,p_8=f_{\widetilde{\tau},8}(\hat{p}_8),p_9=f_9\bigl(\hat{t}\bigr) \big \}\subset \mathbb{C}/ \langle 1, \widetilde{\tau} \rangle
 \]
 for each ${\hat{t}=(\widetilde{\tau},\hat{p}_1,\dots ,\hat{p}_8)\in T}$, a~11-dimensional complex manifold $\mathcal{F}$ can be derived from $\operatorname{CP}^2 \times T$. Moreover, let $\widetilde{\pi}$ be defined as a map from $\mathcal{F}$ to $T$ such that \smash{$(\widetilde{\pi})^{-1}\bigl(\hat{t}\bigr)$} is the blow up of \smash{$\pi^{-1}\bigl(\hat{t}\bigr)$} at the nine points in~$ { \{ f_{\widetilde{\tau}}(p_1),\dots ,f_{\widetilde{\tau}}(p_9) \} \times \{ (\widetilde{\tau},\hat{p}_1,\dots ,\hat{p}_8) \}}$ with
\[
\{ p_1=f_{\widetilde{\tau},1}(\hat{p}_1),\dots , p_8=f_{\widetilde{\tau},8}(\hat{p}_8),p_9=f_9\bigl(\hat{t}\bigr) \}\subset \mathbb{C}/ \langle 1, \widetilde{\tau} \rangle
\]
for each $\hat{t}=(\widetilde{\tau},\hat{p}_1,\dots ,\hat{p}_8)\in T$. Now the normal bundle of the strict transform of the elliptic curve $f_{\widetilde{\tau}}(C_0(\widetilde{\tau}))\times \{ (\widetilde{\tau}, \hat{p}_1,\dots ,\hat{p}_8) \}$ in $(\widetilde{\pi})^{-1}\bigl(\hat{t}\bigr)$ satisfies the Diophantine condition corresponding to $(p,q)$ for $\hat{t}=(\widetilde{\tau},\hat{p}_1,\dots ,\hat{p}_8)\in T$ \cite{Takayuki-Takato2}. Let $\mathcal{S} \subset \mathcal{F}$ be the strict transform of~${\mathcal{S}_0 \subset \operatorname{CP}^2 \times T}$. Then the following lemma is true.

\begin{Lemma}[{\cite{Takayuki-Takato1}}] \label{Lemma3.3}
 Through blowing up nine points on every fiber of \emph{$\bigl(\operatorname{CP}^2 \times T,T, \pi\bigr)$} such that the normal bundle of the strict transform of the elliptic curve mentioned above embedded in each fiber satisfies the Diophantine condition corresponding to the Diophantine number pair $(p,q)$, the complex analytic families \emph{$(\mathcal{F},T, \widetilde{\pi})$} and \emph{$(\mathcal{S},T, \widetilde{\pi}|_\mathcal{S})$} can be constructed.
\end{Lemma}

The complete proof of Lemma~\ref{Lemma3.3} is in Appendix~\ref{appendixA}. Here, $(\mathcal{F},T, \widetilde{\pi})$ is a complex analytic family of \smash{$\operatorname{CP}^2 \#9 \overline{\operatorname{CP}}^2$} over the 9-dimensional complex manifold $T$ and $(\mathcal{S},T, \widetilde{\pi}|_\mathcal{S})$ is a complex analytic family of elliptic curves over the 9-dimensional complex manifold $T$.

Secondly, the following is to give a deformation family of quasi-projective varieties.

Here, the fibers of deformation families defined in Definition~\ref{Definition3.1} should be compact. Since there will be a deformation family of open complex manifolds to be confirmed, the definitions of deformation in \cite{Elizabeth-Francisco} will be used. Then the corresponding definition for the category of smooth manifolds is as follows.

\begin{Definition}[the category of smooth manifolds] \label{Definition3.4}
 The objects for the category of smooth manifolds are the smooth manifolds. For any objects $M$ and $N$, the morphisms from $M$ to $N$ are the maps in the set $C^\infty (M,N)$.
\end{Definition}

So the definition of deformation family from \cite{Elizabeth-Francisco} is as follows.

\begin{Definition}[{\cite[Definition~2.3]{Elizabeth-Francisco}}] \label{Definition3.5}
 A deformation family of a complex manifold $X_1$ is a~holomorphic surjective submersion $\mathcal{X} \rightarrow D_1$, where $D_1$ is a complex disc centred at $0$, satisfying
\begin{itemize}\itemsep=0pt
\item[$(1)$] $\pi^{-1}(0)=X_1$,
\item[$(2)$] $\mathcal{X}$ is locally trivial in the $C^\infty$ category.
\end{itemize}
For any $t\in D_1$, the fiber $\pi^{-1}(t)$ is called a deformation of $X_1$.
\end{Definition}

The discussion about Definition~\ref{Definition3.5} can be found in \cite{Edoardo-Elizabeth-Francisco}. Here, there is a simple corresponding lemma.

\begin{Lemma} \label{Lemma3.6}
For a complex analytic family $(\mathcal{M}_1, B_1, \varpi)$ with $B_1$ being a complex disc centred at $0$, $(\mathcal{M}_1, B_1, \varpi)$ is also a deformation family of the complex manifold $\varpi^{-1}(0)$ satisfying the conditions in Definition {\rm\ref{Definition3.5}}.
\end{Lemma}

\begin{proof}
From the Definition~\ref{Definition3.1}, since the rank of the Jacobian matrix of $\varpi$ is equal to 1 at every point of $\mathcal{M}_1$, $\varpi$ is a holomorphic surjective submersion $\mathcal{M}_1 \rightarrow B_1$. Moreover, $\varpi^{-1}(0)$ is a compact complex submanifold of $\mathcal{M}_1$. In addition, the following lemma is from \cite{Kunihiko}.

\begin{Lemma}[{\cite[Theorem~2.5]{Kunihiko}}] \label{Lemma3.7} $\forall x \in B_1$, $\exists$ a polydisc $U_x \subset B_1$ centered at $x$ with a diffeomorphism $\Psi_x \colon \varpi^{-1}(x) \times U_x \rightarrow \varpi^{-1}(U_x)$ such that $\varpi \circ \Psi_x \colon \varpi^{-1}(x) \times U_x \rightarrow U_x$ is a projection.
\end{Lemma}

Taking use of Lemma~\ref{Lemma3.7}, $\varpi$ is $C^\infty$ locally trivial on $B_1$. That is to say, $\mathcal{M}_1$ is locally trivial in the $C^\infty$ category. Therefore, $(\mathcal{M}_1, B_1, \varpi)$ is a deformation family of the complex manifold~$\varpi^{-1}(0)$ satisfying the conditions in Definition~\ref{Definition3.5}. So for a complex analytic family~$(\mathcal{M}_1, B_1, \varpi)$ with $B_1$ being a complex disc centred at $0$, $(\mathcal{M}_1, B_1, \varpi)$ is also a deformation family of the complex manifold $\varpi^{-1}(0)$ satisfying the conditions in Definition~\ref{Definition3.5}. Then Lemma~\ref{Lemma3.6} is proved.
\end{proof}

In addition, the following is to introduce the definition of the smooth compacitifiable deformation family. The complex model space is defined as follows.

\begin{Definition}[{\cite[Definition~13]{Edoardo-Elizabeth-Francisco}}] \label{Definition3.8}
$\forall f_1,\dots ,f_k \in \mathcal{O}(D)$ with $D \subset \mathbb{C}^n$ being a domain and~${n \in \mathbb{N}}$, then $I_D=\mathcal{O}_Df_1+\mathcal{O}_Df_2+\dots +\mathcal{O}_Df_k$ is an ideal sheaf. Let $X_2:=N(f_1,\dots ,f_k)$ and $\mathcal{O}_{X_2}=(\mathcal{O}_D/I_D)|_{X_2}$. So the $\mathbb{C}$-ringed space $(X_2,\mathcal{O}_{X_2})$ is called a complex model space.
\end{Definition}

The complex space is defined as follows.

\begin{Definition}[{\cite[Definition~14]{Edoardo-Elizabeth-Francisco}}] \label{Definition3.9}
For a $\mathbb{C}$-ringed space $\bigl(\hat{X},\mathcal{O}_{\hat{X}}\bigr)$, if $\hat{X}$ is a Hausdorff space and \smash{$\forall x \in \hat{X}$}, $\exists$ an open neighborhood $U_x$ such that the open $\mathbb{C}$-ringed subspace $(U_x,\mathcal{O}_{U_x})$ is isomorphic to a complex model space, then $\bigl(\hat{X},\mathcal{O}_{\hat{X}}\bigr)$ is called a complex space.
\end{Definition}

So the following is the definition of the smooth compactifiable deformation family.

\begin{Definition}[{\cite[Definition~16]{Edoardo-Elizabeth-Francisco}}] \label{Definition3.10}
For a complex space $X$, $X$ is compactifiable if there is an open embedding $\rho\colon X\rightarrow \widetilde{X}$ such that $\rho(X)=\widetilde{X} \backslash D_2$ with $\widetilde{X}$ being a compact complex space and $D_2$ being a closed analytic subset of $\widetilde{X}$. Then a compactifiable deformation of $X$ is a~sixtuple~$(\mathcal{X},D_2,Y_1,f,o,j)$ such that
\begin{itemize}\itemsep=0pt
\item[$(1)$] $\mathcal{X}$ and $Y_1$ are complex spaces,
\item[$(2)$] $D_2$ is a closed analytic subset of $\mathcal{X}$ and $o$ is a point in $Y_1$,
\item[$(3)$] $f$ is a proper holomorphic map from $\mathcal{X}$ to $Y_1$ such that $f|_{\mathcal{X} \backslash D_2}$ is flat,
\item[$(4)$] $j$ is an open embedding from $X$ to $f^{-1}(o)$ with $j(X)=f^{-1}(o)\backslash D_2$.
\end{itemize}
\end{Definition}

Furthermore, the differentially local triviality can be defined as follows.

\begin{Definition}[{\cite[Definition~21]{Edoardo-Elizabeth-Francisco}}] \label{Definition3.11}
For a holomorphic submersion $g\colon X \rightarrow \hat{Y}$ with $X$ and~$\hat{Y}$ being complex manifolds, $g$ is differentially locally trivial over $\hat{Y}$ if $\forall y \in \hat{Y}$, $\exists$ an open neighborhood $U_y$ such that $g|_{g^{-1}(U_y)}$ is differentially isomorphic to the projection ${g^{-1}(y) \times U_y \rightarrow U_y}$.
\end{Definition}

Then for the compactifiable deformation of $X$ defined above, if $f|_{\mathcal{X} \backslash D_2 }\colon \mathcal{X} \backslash D_2\rightarrow Y_1$ with $\mathcal{X} \backslash D_2$ and $Y_1$ being complex manifolds is differentially locally trivial over $Y_1$, then $(\mathcal{X}, D_2, Y_1,f,\allowbreak o,j)$ is said to be differentially trivial along $Y_1$.

\begin{Definition} \label{Definition3.12}
For a compactifiable complex space $X$ defined above, assuming that $X$ is smooth, if $f|_{\mathcal{X}\backslash D_2}$ is a smooth morphism of complex spaces and each fiber of $f$ is a complex manifold, then $(\mathcal{X}, D_2, Y_1,f, o, j)$ is a smooth compactifiable deformation.
\end{Definition}

Taking use of the notations above, let $t_0 = (\tau_0, q_1,q_2,\dots ,q_8) \in T$, then
\[\{ p_1=f_{\tau_0,1}(q_1),\dots ,p_8=f_{\tau_0,8}(q_8),p_9=f_9(t_0) \} \subset C_0(\tau_0)=C_0.\]
Let \smash{$S \cong \operatorname{CP}^2 \#9\overline{\operatorname{CP}}^2$} be the blow up of $S_0 =\operatorname{CP}^2$ at the nine points in the set
\[ \{ f_{\tau_0}(p_1),\dots ,f_{\tau_0}(p_9) \}\subset f_{\tau_0}(C_0)\]
and $C$ be the strict transform of $f_{\tau_0}(C_0)$. Here, the normal bundle $N_{C/S}$ should satisfy the Diophantine condition corresponding to the Diophantine pair $(p,q)$.

A smooth compactifiable deformation of the quasi-projective variety \smash{$S \backslash C \subset S \cong \operatorname{CP}^2 \#9\overline{\operatorname{CP}}^2$} differentially trivial along $T$ is as follows.

\begin{Theorem} \label{Theorem3.13}
Let $o=(\tau_0,q_1,\dots ,q_8)$. Define $j\colon S \backslash C\rightarrow \widetilde{\pi}^{-1}(o)$ by $j(x)=(x,\tau_0, q_1,q_2,\dots ,q_8)$ for $x\in S\backslash C$. Then $(\mathcal{F},\mathcal{S},T,\widetilde{\pi},o,j)$ is a compactifiable deformation of $S\backslash C$.
\end{Theorem}

\begin{proof}
Firstly, it is to prove that $S\backslash C$ is compactifiable.
$S\backslash C$ is a quasi-projective variety. $(\widetilde{\pi})^{-1}(o)\cong S$ is a compact complex manifold. So $(\widetilde{\pi})^{-1}(o)$ is a complex space. Moreover, $(\widetilde{\pi})^{-1}(o)\cap \mathcal{S}\cong \mathbb{C}/\langle 1, \tau_0 \rangle$ is a closed analytic submanifold of $(\widetilde{\pi})^{-1}(o)$. Since $(\widetilde{\pi})^{-1}(o)\backslash j(S\backslash C)=(\widetilde{\pi})^{-1}(o)\cap \mathcal{S}$ is closed, $j(S\backslash C)$ is an open submanifold of~$(\widetilde{\pi})^{-1}(o)$. Then $j$ is an open embedding with $j(S\backslash C)=(\widetilde{\pi})^{-1}(o)\backslash \mathcal{S}$. So $((\widetilde{\pi})^{-1}(o),\widetilde{\pi}^{-1}(o)\cap \mathcal{S},j)$ is a compactification of $S\backslash C$. That is to say, $S\backslash C$ is compactifiable.

Secondly, it is to prove that $\mathcal{S}$ is a closed analytic subset of $\mathcal{F}$ and $\widetilde{\pi}$ is a proper holomorphic map from $\mathcal{F}$ to $T$ such that \smash{$\widetilde{\pi}|_{\mathcal{F}\backslash \mathcal{S}}$} is flat.

Since $\mathcal{F}$ and $T$ are complex manifolds, $\mathcal{F}$ and $T$ are complex spaces. From above, $\mathcal{S}$ is a~codimension 1 complex submanifold of $\mathcal{F}$. In addition, since the blow up points are all on $\mathcal{S}_0$, it is easy to get an open neighborhood for every point in $\mathcal{F}\backslash \mathcal{S}$. Therefore, $\mathcal{F}\backslash \mathcal{S}$ should be an open subset of $\mathcal{F}$. Then $\mathcal{S}$ is a closed subset of $\mathcal{F}$. So $\mathcal{S}$ is a closed analytic subset of $\mathcal{F}$.

Furthermore, since $(\mathcal{F},T,\widetilde{\pi})$ is the complex analytic family of \smash{$\operatorname{CP}^2 \#9\overline{\operatorname{CP}}^2$}, $\widetilde{\pi}$ is a proper holomorphic map from $\mathcal{F}$ to $T$. So $\widetilde{\pi}|_{\mathcal{F}\backslash \mathcal{S}}$ is also proper.

Since $\widetilde{\pi}|_{\mathcal{F}\backslash \mathcal{S}}$ is a proper submersion from the complex manifold $\mathcal{F}\backslash \mathcal{S}$ to $T$, $\widetilde{\pi}|_{\mathcal{F}\backslash \mathcal{S}}$ is flat.

Finally, it is to prove that $(\mathcal{F},\mathcal{S},T,\widetilde{\pi},o,j)$ is a compactifiable deformation of $S\backslash C$.

From above, $j\colon S\backslash C \rightarrow (\widetilde{\pi})^{-1}(o)$ is an open embedding with $j(S\backslash C)=(\widetilde{\pi})^{-1}(o)\backslash \mathcal{S}$. Therefore, with all the proofs above, $(\mathcal{F},\mathcal{S},T,\widetilde{\pi},o,j)$ is a compactifiable deformation of $S\backslash C$. So Theorem~\ref{Theorem3.13} is proved.
\end{proof}

In addition, the following statement is true.

\begin{Proposition} \label{Proposition3.14}
$(\mathcal{F},\mathcal{S},T,\widetilde{\pi},o,j)$ is a smooth compactifiable deformation differentially trivial along $T$.
\end{Proposition}

\begin{proof}
Since \smash{$\widetilde{\pi}|_{\mathcal{F}\backslash \mathcal{S}}$} is a proper submersion from the complex manifold $\mathcal{F}\backslash \mathcal{S}$ to the complex manifold $T$ and \smash{$(\widetilde{\pi})^{-1}\bigl(\hat{t}\bigr)\cong \operatorname{CP}^2 \#9\overline{\operatorname{CP}}^2$} is a compact complex manifold for $\hat{t} \in T$, $(\mathcal{F},\mathcal{S},T,\widetilde{\pi},o,j)$ is a smooth compactifiable deformation.

Moreover, the following lemma from Ehresmann indicates that $(\mathcal{F},\mathcal{S},T,\widetilde{\pi},o,j)$ is differentially trivial along $T$.

\begin{Lemma}[{\cite[Theorem~9.3]{Claire}}] \label{Lemma3.15}
 For any proper holomorphic submersion $g\colon M\rightarrow N$ with $M$ and $N$ being complex manifolds, $g$ is a locally trivial fibration.
\end{Lemma}

Taking use of Lemma~\ref{Lemma3.15}, since $\widetilde{\pi}|_{\mathcal{F}\backslash \mathcal{S}}$ is a proper submersion from the complex manifold~$\mathcal{F}\backslash \mathcal{S}$ to the complex manifold $T$, $\widetilde{\pi}|_{\mathcal{F}\backslash \mathcal{S}}$ is a locally trivial fibration. So $(\mathcal{F},\mathcal{S},T,\widetilde{\pi},o,j)$ is differentially trivial along $T$.

In conclusion, $(\mathcal{F},\mathcal{S},T,\widetilde{\pi},o,j)$ is a smooth compactifiable deformation differentially trivial along $T$. So Proposition \ref{Proposition3.14} is proved.
\end{proof}

Therefore, $(\mathcal{F},\mathcal{S},T,\widetilde{\pi},o,j)$ is a compactifiable deformation of the quasi-projective variety $S\backslash C$ over a 9-dimensional complex manifold $T$ with $o=(\tau_0,q_1,q_2,\dots ,q_8)$ and $j\colon S \backslash C \rightarrow (\widetilde{\pi}^{-1}(o))$ defined by $j(x)=(x,\tau_0,q_1,q_2,\dots ,q_8)$ for $x \in S\backslash C$. Combining the results of Theorem~\ref{Theorem3.13} and Proposition \ref{Proposition3.14}, the proof of Theorem~\ref{Theorem1.1} is completed.

\subsection[Complete K\"ahler metrics on each fiber of the deformation family of quasi-projective varieties]{Complete K\"ahler metrics on each fiber of the deformation family\\ of quasi-projective varieties}

The following is to construct the complete K\"ahler metrics on each fiber of the deformation family of quasi-projective varieties. A key point for the construction here is to take use of Macro Brunella's statement in \cite{Marco}.

Taking use of the notations from the introduction and Section~\ref{sec2}, the complete K\"ahler metrics will be constructed on $S(\tau)\backslash C(\tau)$ for $\tau \in Y$.

For any integers $k\geq 2$ and $d \geq 3k+1$ satisfying $10 \leq d^2/k^2$, let $L=d\cdot H-k\cdot\sum _{i=1}^9 E_i$. Taking use of Lemma~\ref{Lemma2.5}, $L$ is an ample line bundle on $S(\tau)$. Then there exists $n\in \mathbb{N}$ such that~$L^n$ and $L^n\otimes [-C(\tau)]$ are very ample \cite{Phillip}. Now let $ \{ g_0,g_1,\dots ,g_{\theta} \}$ be a basis for the space of holomorphic sections on $L^n\otimes [-C(\tau)]$ with $\theta \in \mathbb{N}$. Let $\hat{s}$ be the defining section of $C(\tau)$ and~$ \{ \zeta \}$ be a frame for $L$. Then $ \{ g_0\hat{s},g_1\hat{s},\dots ,g_{\theta} \hat{s} \}$ are all holomorphic sections on $L^n$ with value zero on~$C(\tau)$. This is an important point for the construction. Then the hermitian metric $h_L$ defined on $L|_{S(\tau)\backslash C(\tau)}$ through the $C^\infty $ function
\[
h_L(x)=\frac{\zeta(x) \cdot \overline{\zeta (x)}}{\bigl(|g_0(x)\cdot \hat{s}(x)|^2+|g_1(x)\cdot \hat{s}(x)|^2+\dots +|g_{\theta}(x)\cdot \hat{s}(x)|^2\bigr)^{\frac{1}{n}}}
\]
for $x \in S(\tau)\backslash C(\tau)$ is smooth.

\begin{Lemma} \label{Lemma3.16}
The curvature form $\Theta_{h_L}$ of $h_L$ is positive.
\end{Lemma}

\begin{proof}
Since $L^n \otimes [-C(\tau)]$ is very ample, taking use of the Kodaira embedding theorem \cite{Phillip}, the evaluation map
\[\Phi_{L^n \otimes [-C(\tau)]}(x)=[g_0(x):g_1(x):\dots :g_{\theta}(x)]\in \operatorname{CP}^{\theta}\]
gives a holomorphic embedding of \smash{$S(\tau) \cong \operatorname{CP}^2 \#9\overline{\operatorname{CP}}^2$} into $\operatorname{CP}^{\theta}$. In addition, let $[\xi_0,\xi_1,\dots ,\xi_{\theta}]$ be the coordinate system of $\operatorname{CP}^{\theta}$ coinciding with the image of the evaluation map $\Phi_{L^n \otimes [-C(\tau)]}$.

The functions
\[
g_i([\xi_0:\xi_1:\dots :\xi_{\theta}])=\frac{|\xi_i|^2}{|\xi_0|^2+|\xi_1|^2+\dots +|\xi_{\theta}|^2}
\]
on $\big \{ [\xi_0:\xi_1:\dots :\xi_{\theta}]\in \operatorname{CP}^\theta\mid \xi_i\neq 0\big \}$ for $i \in \{ 0,1,\dots ,\theta \}$ define a hermitian metric $g$ on the hyperplane line bundle $\mathcal{O}(1)$. Then the curvature form for the hermitian metric $g$ is
\[
\Theta_g=\frac{\sqrt{-1}}{\pi}\partial\overline{\partial}\log\bigl(|\xi_0|^2+|\xi_1|^2+\dots +|\xi_{\theta}|^2\bigr)
\]
for $[\xi_0:\xi_1:\dots :\xi_{\theta}]\in \operatorname{CP}^\theta$. Furthermore, the associate $(1,1)$-form for the Fubini--Study metric on $\operatorname{CP}^\theta$ is the positive K\"ahler form
\[\frac{\sqrt{-1}}{2\pi}\partial\overline{\partial}\log\bigl(|\xi_0|^2+|\xi_1|^2+\dots +|\xi_{\theta}|^2\bigr)\]
for $[\xi_0:\xi_1:\dots :\xi_{\theta}]\in \operatorname{CP}^\theta$. So $\Theta_g$ is a positive form on $\operatorname{CP}^\theta$.

Therefore, $g$ is a smooth positive hermitian metric on the hyperplane line bundle $\mathcal{O}(1)$ for the coordinate system $[\xi_0:\xi_1:\dots :\xi_{\theta}]$ on $\operatorname{CP}^{\theta}$. Then the pullback metric $\hat{g}=(\Phi_{L^n \otimes [-C(\tau)]})^*g$ is a smooth positive hermitian metric on the pullback bundle $(\Phi_{L^n \otimes [-C(\tau)]})^*\mathcal{O}(1)=L^n \otimes [-C(\tau)]$. So the curvature form $\Theta_{\hat{g}}$ of $\hat{g}$ is positive.

Moreover, in \cite{Marco}, Macro Brunella stated that under the condition that the normal bundle~\smash{$N_{C(\tau)/S(\tau)}$} satisfied the Diophantine condition, \smash{$K^{-1}_{S(\tau)}=[C(\tau)]$} admitted a smooth hermitian metric with semi-positive curvature. So $\log|\hat{s}|$ gives the local weights of a (singular) metric $g_{\hat{s}}$ on $K^{-1}_{S(\tau)}$ with $\Theta_{g_{\hat{s}}}>0$ \cite[Definition~8.7]{Vincent-Ahmed}. Therefore, $\Theta_{h_L}=\frac{1}{n}(\Theta_{\hat{g}}|_{S(\tau)\backslash C(\tau)}+\Theta_{g_{\hat{s}}}|_{S(\tau)\backslash C(\tau)})$ is positive. Then Lemma~\ref{Lemma3.16} is proved. \end{proof}

In addition, taking use of Theorem~\ref{Theorem2.19}, let $L|_{C(\tau)}$ be the restricted line bundle of $L$ on $C(\tau)$. Furthermore, let $L^*$ be the pull back bundle of $L|_{C(\tau)}$ on $V_{0,\infty}$. Then $L^*|_{V_{0,r}}=L|_{V_{0,r}}$. Taking use of the equivalent condition in \cite[Theorem]{Christian} (also see \cite{Takayuki-Takato2}), there exists a corresponding theta line bundle on $U_0=\mathbb{C}^2/{\Lambda_0}$ holomorphically isomorphic to the pull back line bundle $(f_{U_0})^*(L^*)$.

Let \smash{$H_1=\bigl(\begin{smallmatrix}
\frac{(L.C(\tau))}{\operatorname{ Im}\tau} & 0\\
0&0
\end{smallmatrix}\bigr)$}. Here, $H_1(x,y)=x^tH_1 \overline{y}$ for $x,y \in \mathbb{C}^2$. Defining
$\rho\colon \Lambda_0 \rightarrow U(1)$ through
\[
\rho(\lambda+\mu)=\rho(\lambda)\rho(\mu){\rm e}^{\pi\sqrt{-1}\cdot\operatorname{ Im}H_1(\lambda,\mu)}
\]
 for $ \lambda, \mu \in \Lambda_0 $, let \[\alpha_\lambda(x)=\rho (\lambda){\rm e}^{\pi H_1(x,\lambda)+(\pi/2)\cdot H_1(\lambda,\lambda)}\] for $ \lambda \in \Lambda_0 $ and $ x \in \mathbb{C}^2 $.

The theta line bundle on $U_0=\mathbb{C}^2/{\Lambda_0}$ corresponding to \smash{$H_1=\bigl(\begin{smallmatrix}
\frac{(L.C(\tau))}{\operatorname{ Im} \tau} & 0\\
0&0
\end{smallmatrix}\bigr)$} can be defined~as ${
L_{H_1,\rho}=\bigl(\mathbb{C}_{\zeta} \times \mathbb{C}^2\bigr)/{\Lambda_0}}
$
with $\lambda \cdot (\zeta, x)=$$(\alpha_\lambda(x)\cdot \zeta, x+\lambda )$ for $\lambda \in\Lambda_0$, $\zeta \in \mathbb{C}_\zeta$ and $x \in \mathbb{C}^2$.

\begin{Lemma} \label{Lemma3.17}
The theta line bundle $L_{H_1,\rho}$ on $U_0=\mathbb{C}^2/{\Lambda_0}$ is holomorphically isomorphic to the pull back bundle $(f_{U_0})^*(L^*)$ with \smash{$H_1=\bigl(\begin{smallmatrix}
\frac{(L.C(\tau))}{\operatorname{ Im} \tau} & 0\\
0&0
\end{smallmatrix}\bigr)$}.
\end{Lemma}

\begin{proof}
For the pull back bundle \smash{$L^*_{H_1,\rho}=\bigl(f^{-1}_{U_0}\bigr)^*L_{H_1,\rho}$} on $V_{0,\infty}$, according to \cite[Proposition~3.5]{Takayuki-Takato2}, $\exists \hat{G} \in \operatorname{Pic}(W) $ such that \smash{$\hat{G}|_{V_{0,r}}=L^*_{H_{1,\rho}}|_{V_{0,r}}$}. Then \smash{$(\hat{G}.C(\tau))=(L.C(\tau))$} (see \cite[Lemma~3.4]{Takayuki-Takato2}).

Moreover, according to \cite[Lemma~3.4]{Takayuki-Takato2} and \cite[Proposition~3.5]{Takayuki-Takato2}, $L_{H_1,\rho}$ with $H_1=\smash{\bigl(\begin{smallmatrix}
\frac{(L.C(\tau))}{\operatorname{ Im}\tau} & 0\\
0&0
\end{smallmatrix}\bigr)}$ is the unique theta line bundle defined in the formula as above such that $\hat{G}$ has the same intersection number with $C(\tau)$ as $L$. Therefore, the theta line bundle $L_{H_1, \rho}$ on $U_0=\mathbb{C}^2_{(z,\eta)}/{\Lambda_0}$ is holomorphically isomorphic to the pull back bundle $(f_{U_0})^*(L^*)$ with \smash{$H_1\!=\!\bigl(\!\begin{smallmatrix}
\frac{(L.C(\tau))}{\operatorname{ Im}\tau} & 0\\
0&0
\end{smallmatrix}\bigr)$}. So Lemma~\ref{Lemma3.17} is proved.
\end{proof}

Taking use of Lemma~\ref{Lemma3.17},
\[
c_1(L|_{C(\tau)})=\left[\frac{(L.C(\tau))}{\operatorname{Im}\tau}\cdot \sqrt{-1}\text{{\rm d}z}\land {\rm d}\bar{z}\right].
\]
In addition, according to \cite[Lemma~7.31]{Vincent-Ahmed} ($\partial \overline{\partial}$-lemma) and \cite[Section 8.2]{Vincent-Ahmed}, there exists a~smooth hermitian metric $h_{C(\tau)}$ on $L|_{C(\tau)}$ such that
\[
\Theta_{h_{C(\tau)}} = \frac{(L.C(\tau))}{\operatorname{ Im}\tau} \sqrt{-1} {\rm d}z \wedge {\rm d} \bar{z}.
\]
 So $(\pi_W)^*h_{C(\tau)}$ is a smooth hermitian metric on $L|_W$ such that
\[
\Theta_{(\pi_W)^*h_{C(\tau)}}=\frac{(L.C(\tau))}{\operatorname{Im}\tau} \sqrt{-1} {\rm d}z \wedge d \bar{z}.
\]

For the line bundle \smash{$L|_{C(\tau)}$}, there exist open finite covers \smash{$\{ \mathcal{U}_{\alpha_1} \}_{\alpha_1 \in I}$} of $C(\tau)$ and local trivializations
\[
\Phi_{\alpha_1}\colon \ \pi^{-1}_{L|_{C(\tau)}}(\mathcal{U}_{\alpha_1}) \rightarrow \mathcal{U}_{\alpha_1} \times \mathbb{C}
\]
 mapping \smash{$\pi^{-1}_{L|_{C(\tau)}}(x)$} isomorphically onto $ \{ x \} \times \mathbb{C}$ for $x \in \mathcal{U}_{\alpha_1}$ and $\alpha_1 \in I$ with $I$ being an index set. Let
\[
g_{\alpha_1 \beta_1}(x)=\big[\Phi_{\alpha_1} \circ (\Phi_{\beta_1})^{-1}\big]|_{ \{ x \} \times \mathbb{C}}
\]
for $ x \in \mathcal{U}_{\alpha_1} \cap \mathcal{U}_{\beta_1}$ with $\alpha_1, \beta_1 \in I$.

The line bundle \smash{$L|_W=(\pi_W)^*L|_{C(\tau)}$} and \smash{$\big \{ W_{\alpha_1}=\pi^{-1}_W(\mathcal{U}_{\alpha_1})\big \}_{\alpha_1 \in I}$} is an open cover of $W$. Moreover, there exist local trivializations
\[
\Psi_{\alpha_1}\colon \ \pi^{-1}_{L|_W}(W_{\alpha_1})\rightarrow W_{\alpha_1} \times \mathbb{C}
\]
 mapping \smash{$\pi^{-1}_{L|_W}(x)$} isomorphically onto $\{ x \} \times \mathbb{C}$ induced by $\Phi_{\alpha_1}$ for $x=[(z,w)] \in W_{\alpha_1}$ and $\alpha_1 \in I$. Let
\[
\widetilde{g}_{\alpha_1 \beta_1}(x)=\big[\Psi_{\alpha_1} \circ (\Psi_{\beta_1})^{-1}\big]\big|_{\{ x \}\times \mathbb{C}}
\]
for $x \in W_{\alpha_1} \cap W_{\beta_1}$ with $\alpha_1,\beta_1 \in I$.

A singular hermitian metric for a line bundle is defined through local weights \cite[Definition~8.7]{Vincent-Ahmed}. For $\varepsilon > 0$, let
$
 \varphi^{\alpha_1}_L \in \operatorname{Psh}(W_{\alpha_1} \backslash \mathcal{U}_{\alpha_1})\cap C^\infty(W_{\alpha_1} \backslash \mathcal{U}_{\alpha_1})$ and $\varphi^{\alpha_1} _{C(\tau)}+\log \varepsilon \in \operatorname{Psh}(W_{\alpha_1}) \cap C^\infty (W_{\alpha_1})
$
be the local weights of the metrics $h_L$ and \smash{$\varepsilon^{-1}(\pi_W)^*h_{C(\tau)}$} with $\operatorname{Psh}(W_{\alpha_1} \backslash \mathcal{U}_{\alpha_1})$ and $\operatorname{Psh}(W_{\alpha_1})$ being the sets consisting of plurisubharmonic functions on $W_{\alpha_1} \backslash \mathcal{U}_{\alpha_1}$ and $W_{\alpha_1}$ for $\alpha_1 \in I$, respectively.

Here, it is another important point for the construction modified from \cite{Takayuki-Takato2}. Through the following theorem, $\varphi^{\alpha_1}_L$ and $\varphi^{\alpha_1}_{C(\tau)}+\log \varepsilon$ can be patched to smooth plurisubharmonic functions on $W_{\alpha_1}$ for $\alpha_1 \in I$.

Let $\gamma =(\gamma_1,\gamma_2,\dots ,\gamma_p)$ with $\gamma_i >0$ for $i \in \left \{ 1,\dots ,p\right \}$ and $X$ be a compact complex manifold. Moreover, let $\hat{\theta}$ be a real non-negative smooth function defined on $ \mathbb{R}$ supported on $[-1,1]$ with
\[
\int _\mathbb{R} \hat{\theta}(x){\rm d}x=1 \qquad \text{and} \qquad \int _\mathbb{R}x\hat{\theta} (x){\rm d}x=0.
\]
 In addition, let
\[
M_\gamma (t_1,t_2,\dots ,t_p)=\int _{\mathbb{R}^p}\max \{ t_1+h_1,\dots ,t_p+h_p\}\prod_{1 \leq j \leq p} \gamma^{-1}_j \hat{\theta} (h_j/\gamma_j) {\rm d}h_1 {\rm d}h_2 \dots .{\rm d} h_p
\]
for $(t_1,t_2,\dots ,t_p) \in \mathbb{R}^p$ which is called a regularized max function \cite{Alexandre}.

\begin{Lemma}[{\cite[Lemma~5.18]{Jean}}] \label{Lemma3.18}
 $M_\gamma (t_1,t_2,\dots ,t_p)$ defined above satisfying the following properties:
\begin{itemize}\itemsep=0pt
\item[$(1)$]$M_\gamma (t_1,t_2,\dots ,t_p)$ is smooth and convex in $\mathbb{R}^p$;
\item[$(2)$]$M_\gamma (t_1+a,t_2+a,\dots ,t_p+a)=M_\gamma (t_1,t_2,\dots ,t_p)+a$, for any $a \in \mathbb{R}$ and $(t_1,t_2,\dots ,t_p) \in \mathbb{R}^p$.
\end{itemize}
\end{Lemma}

Let $h=h_L$ on $S(\tau)\backslash W$. Moreover, let
\[
M_{(1,1)}\bigl(\varphi^{\alpha_1}_L,\varphi^{\alpha_1}_{C(\tau)}+\log\varepsilon\bigr)=\int_{\mathbb{R}^2}\max\bigl\{\varphi^{\alpha_1}_L+h_1,\varphi^{\alpha_1}_{C(\tau)}
+\log\varepsilon+h_2 \bigr\}\prod_{1\leq j\leq 2}\hat{\theta}(h_j){\rm d}h_1{\rm d}h_2\]
be the local weight of the metric $h$ on $W_{\alpha_1}$ for $\alpha_1 \in I$.

Let $r_1$ and $r_2$ be positive numbers with $r_1 <r_2<r$. Choosing $\varepsilon$ small enough,
\begin{gather*}
M_{(1,1)}\bigl(\varphi^{\alpha_1}_L,\varphi^{\alpha_1}_{C(\tau)}+\log \varepsilon\bigr)\\
\qquad=\int_{\mathbb{R}^2}\max \bigl\{\varphi^{\alpha_1}_L+h_1,\varphi^\alpha_{C(\tau)}+\log\varepsilon +h_2 \bigr\} \prod_{1 \leq j \leq 2} \hat{\theta} (h_j){\rm d} h_1 {\rm d} h_2 =\varphi^{\alpha_1}_L
\end{gather*}
on $W_{\alpha_1} \cap \{ [(z,w)] \in W \mid r_1 <|w| \}$ for $\alpha_1 \in I$. That is to say, $h=h_L$ on $ \{ [(z,w)] \in W \mid r_1<|w| \}$.

Taking use of property~(1) of Lemma~\ref{Lemma3.18}, $M_{(1,1)}\bigl(\varphi^\alpha_L,\varphi^\alpha_{C(\tau)}+\log \varepsilon\bigr)$ is plurisubharmonic for~${\alpha\in I}$. Furthermore, from property~(2) of Lemma~\ref{Lemma3.18}, $h$ is a smooth hermitian metric on~$L$. Therefore, $h$ is a smooth hermitian metric on~$L$ with $\Theta_{h}\geq 0 $ (see \cite[Section~8.2]{Vincent-Ahmed} for the definition of singular hermitian metrics).

Furthermore, there exists $\varepsilon_0 >0$ satisfying $\sqrt{\varepsilon_0}r<r_1$ such that
\[
M_{(1,1)}\bigl(\varphi^{\alpha_1}_L,\varphi^{\alpha_1}_{C(\tau)}+\log \varepsilon\bigr)=\varphi^{\alpha_1}_L
\]
on $W_{\alpha_1} \cap \{ [(z,w)] \in W\mid |w|< \sqrt{\varepsilon_0}r \}$ for $\alpha_1 \in I$. So
$h=\varepsilon^{-1}(\pi_W)^*h_{C(\tau)}$ on $\{ [(z,w)] \in W \mid |w|< \sqrt{\varepsilon_0} r \}$.

Let $\theta_s(t)=\bigl(\log\bigl(t^2/|s|\bigr)\bigr)^2$ be a smooth function for $t>0$ and $s \in \{ s\mid s \in \Delta \backslash \{ 0 \},\, |s|<\varepsilon_0 \}$ with $\Delta = \{ s \in \mathbb{C} \mid |s|<1 \}$.

Let $f(x)=1$ for $|x|<r-\frac{r-r_2}{2}$ and $f(x)=0$ for $r-\frac{r-r_2}{2} \leq|x|$. The standard mollifier~${\eta \in C^\infty(\mathbb{R})}$ is defined as
\[\eta(x):=\begin{cases}
a\cdot {\rm e}^{1/(|x|^2-1)} & \text{if}\  |x|<1, \\
0 & \text{if}\  |x| \geq 1,
\end{cases}\]
with \smash{$a=1/\int^{1}_{-1}{\rm e}^{1/(|x|^2-1)}{\rm d}x>0$} \cite{Lawrence}.

Moreover, let $\eta_\delta(x):=\frac{1}{\delta}\eta(\frac{x}{\delta})$. Then for $\delta=\frac{r-r_2}{4}$, let
\[\widetilde{f}(x)=\eta_\delta * f(x)=\int^\infty_{-\infty}\eta_\delta(x-y)f(y){\rm d}y.\]
Then $\widetilde{f}$ is a smooth function on $\mathbb{R}$ with $\widetilde{f}(x)=1$ for $|x|<r_2$ and $\widetilde{f}(x)=0$ for $r<|x|$.

Let $\Psi_s\colon S(\tau) \backslash C(\tau) \rightarrow \mathbb{R}$ be the smooth function defined by
\[\Psi_s(p):= \begin{cases}
\widetilde{\theta}_s(|w|), & p=[(z,w)] \in W\backslash C(\tau), \\
\widetilde{\theta}_s(r), & p \notin W
\end{cases}\]
with $\widetilde{\theta}_s(x)=\widetilde{f}(x) \cdot \theta_s(x)$ for $x>0$.

Then the following is to give the modified complete K\"ahler metrics on $S(\tau)\backslash C(\tau)$ analogous to the complete K\"ahler metrics described in \cite{Takayuki-Takato2}.

\begin{Theorem}[{\cite{Takayuki-Takato2}}] \label{Theorem3.19}
 There exists $b>0$ such that the metric $h \cdot {\rm e}^{-b \Psi_s}$ is a smooth hermitian metric on $L|_{S(\tau)\backslash C(\tau)}$ with $\Theta_{h \cdot {\rm e}^{-b \Psi_s}}>0$. Moreover, the K\"ahler form $\Theta_{h \cdot {\rm e}^{-b \Psi_s}}$ is of the form
\[\omega|_{W_{\varepsilon_0}}=\frac{ (L.C(\tau))}{\operatorname{ Im}\tau}\cdot \sqrt{-1} {\rm d}z \wedge {\rm d} \bar{z}+\frac{2b}{\pi} \cdot \frac{\sqrt{-1}{\rm d}w \wedge {\rm d} \bar{w}}{|w|^2}\]
on the set $W_{\varepsilon_0}=\{ [(z,w)]\in W \backslash C(\tau)\mid |w|<\sqrt{\varepsilon_0}r \}$. In addition, $\omega|_{W_{\varepsilon_0}}$ is Ricci-flat. Then the K\"ahler form $\Theta_{h\cdot {\rm e}^{-b\Psi_s}}$ gives a complete K\"ahler metric on $S(\tau)\backslash C(\tau)$.
\end{Theorem}

The complete proof of Theorem~\ref{Theorem3.19} can be found in Appendix~\ref{appendixA}. Therefore, the quasi-projective variety above tends to be a non-compact complete K\"ahler manifold.

\subsection{A deformation family of non-compact complete K\"ahler manifolds}

The following is to show that the deformation family $(\mathcal{F},\mathcal{S},T,\widetilde{\pi},o,j)$ can be easily changed to a deformation family satisfying the conditions in Definition~\ref{Definition3.5}.

Since the fibers are now proved to be non-compact complete K\"ahler manifolds, the deformation family above is proved to be a deformation family of non-compact complete K\"ahler manifolds over a 9-dimensional complex manifold.
Moreover, for $T$ defined above, there exists a complex disc $D_3 \subset T$ centred at $o=(\tau_0,q_1,q_2,q_3,q_4,q_5,q_6,q_7,q_8)$.

\begin{Corollary} \label{Corollary3.20}
Let $f\colon T \rightarrow \mathbb{C}^9$ be defined by $f(x)=x-o$. Then $f \circ \widetilde{\pi}|_{\mathcal{F}\backslash \mathcal{S}}\colon\mathcal{F} \rightarrow f(D_3)$ is a~deformation family of the non-compact complete K\"ahler manifold $j(S\backslash C)$ over a 9-dimensional complex manifold satisfying the conditions in Definition {\rm\ref{Definition3.5}}.
\end{Corollary}

The complete proof of Corollary~\ref{Corollary3.20} can be found in Appendix~\ref{appendixA}. Then the construction of the deformation family of the non-compact complete K\"ahler manifold is completed.

\section{A deformation family of symmetric projective K3 surfaces}\label{sec4}

\subsection{Symmetric K\"ahler metrics on the symmetric projective K3 surfaces}

The following is to give symmetric K\"ahler metrics on the symmetric projective K3 surfaces constructed in Section~\ref{sec2}.

In Section~\ref{sec3}, the complete K\"ahler metric is constructed on $S(\tau)\backslash C(\tau)$ for $\tau \in Y$. That is to say, the complete K\"ahler metrics are constructed on fibers of the smooth compactifiable deformation of $S\backslash C$.

Taking use of the notations from Section~\ref{sec2}, the symmetric projective K3 surface is constructed in Theorem~\ref{Theorem2.8} through $S^+(\tau)\backslash C^+(\tau)$ and $S^-(\tau)\backslash C^-(\tau)$.

Since $S^+(\tau)\backslash C^+(\tau)$ is just a copy of $S(\tau)\backslash C(\tau)$, let $\omega^+$ be a copy of $\omega$ on $S(\tau)\backslash C(\tau)$ constructed above. Let $L^{+}=d\cdot H^{+}-k \cdot \sum _{i=1}^9 E_i^{+}$ be defined on $S^+(\tau)$ as a copy of $L$ in Section~\ref{sec3}. Then \smash{$L^{-}=d\cdot H^{-}-k\cdot\sum _{i=1}^9 E_i^{-}$} is an ample line bundle on $S^-(\tau)$ and \smash{$(\ell_{S^+(\tau)})^*L^-=L^+$}. In addition, let $h^+$ be defined on $L^+$ as a copy of $h\cdot {\rm e}^{-b\Psi_s}$ on $L$. Then there is a hermitian metric~$h^-$ defined on $L^-$ satisfying \smash{$h^+=(\ell_{S^+(\tau)\backslash C^+(\tau)})^*h^-$} such that \smash{$\omega^-=\Theta_{h^-}=\bigl(\ell^{-1}_{S^+(\tau)\backslash C^+(\tau)}\bigr)^* \omega^+$} is a~K\"ahler form on $S^-(\tau)\backslash C^-(\tau)$.

Identifying $V_s^+$ and $V_s^-$ through $f_s$, $M^-_s$ and $M^+_s$ can be glued to a K3 surface $X_s$ together with a line bundle $L_s$ derived from the line bundles $L^+$ and $L^-$. Taking use of the result from Section~\ref{sec3}, $\exists \ \varepsilon_0 >0$ such that $\omega$ is well defined on $S(\tau)\backslash C(\tau)$ with the results of Theorem~\ref{Theorem3.19} holding for $s \in \{ s\mid s \in \Delta \backslash \{ 0 \},\, |s|<\varepsilon_0 \}$. Then $\omega^+$ and $\omega^-$ can be glued to be a K\"ahler form $\omega_s$ on $X_s$ for $s \in \{ s\mid s \in \Delta \backslash \{ 0 \},\, |s|<\varepsilon_0 \}$. That is to say, $L_s$ is ample for~${s \in \{ s\mid s \in \Delta \backslash \{ 0 \},\, |s|<\varepsilon_0 \}}$. So $X_s$ is projective for $s \in \{ s\mid s \in \Delta \backslash \{ 0 \}, \, |s|<\varepsilon_0 \}$ (also see~\cite[Theorem~1.2]{Takayuki-Takato2}).

\begin{Definition} \label{Definition4.1}
Let $\mathcal{K}$ be a K3 surface with an ample line bundle $L_\mathcal{K}$ and a non-trivial holomorphic involution $f_\mathcal{K}$ satisfying $f^*_{\mathcal{K}}L_{\mathcal{K}}=L_{\mathcal{K}}$. Then if $\widehat{\omega}$ is a K\"ahler metric on $\mathcal{K}$ such that~${f^*_{\mathcal{K}}\widehat{\omega}=\widehat{\omega}}$, then $\widehat{\omega}$ is said to be a symmetric K\"ahler form on $\mathcal{K}$ corresponding to a symmetric K\"ahler metric.
\end{Definition}

With the discussion above, the following lemma is true.

\begin{Lemma} \label{Lemma4.2}
There exists a symmetric K\"ahler metric on the projective symmetric K3 surface~$X_s$ for $s \in \{ s\mid s \in \Delta \backslash \{ 0 \},\, |s|<\varepsilon_0 \}$.
\end{Lemma}

\begin{proof}
For $s \in \{ s\mid s \in \Delta \backslash \{ 0 \},\, |s|<\varepsilon_0 \}$, since $\omega^+=(\ell_{S^+(\tau)\backslash C^+(\tau)})^*\omega^-$, $\omega_s=F^*_{X_s}\omega_s$. Then~$\omega_s$ is the symmetric K\"ahler form needed for $s \in \{ s\mid s \in \Delta \backslash \{ 0 \},\, |s|<\varepsilon_0 \}$. So there exists a~symmetric corresponding K\"ahler metric on the projective symmetric K3 surface $X_s$ for $s \in \{ s\mid s \in \Delta \backslash \{ 0 \},\, |s|<\varepsilon_0 \}$. Then Lemma~\ref{Lemma4.2} is proved.
\end{proof}

\subsection{A deformation family of symmetric projective K3 surfaces}

The following is to give a deformation family of symmetric projective K3 surfaces over a 10-dimensional complex manifold and complete the proof of Corollary~\ref{Corollary1.2}. With the ample line bundles provided above, the proofs can be completed directly through the maps derived.

Let $\hat{\Delta}=\{ s\mid s \in \Delta \backslash \{ 0 \},\, |s|<\varepsilon_0 \}$. Then there are manifolds
\[\mathcal{M}^+=\bigl(S(\tau)\times \hat{\Delta}\bigr)\backslash \bigl\{ ([(z,w)],s)\in W \times \hat{\Delta} \mid |w|\leq \sqrt{|s|}/r \bigr\}=\bigl\{ (x,s) \in M^+_s \times \hat{\Delta}\bigr\}\]
and
$\mathcal{M}^-=\bigl\{(\ell_{S^+(\tau)}(x),s) \mid (x,s) \in M^+_s \times \hat{\Delta}\bigr\}$.
In addition, let
\[\mathcal{V}^+=\bigl\{ ([(z,w)],s)\in W \times \hat{\Delta} \mid \sqrt{|s|}/r <|w| <\sqrt{|s|} r\bigr\}=\bigl\{ (x,s) \in V^+_s \times \hat{\Delta}\bigr\}\]
and
$\mathcal{V}^-=\bigl\{ \bigl(\ell_{S^+(\tau)}(x),s\bigr) \mid (x,s) \in V^+_s \times \hat{\Delta}\bigr\}$
be the submanifolds of $\mathcal{M}^+$ and $\mathcal{M}^-$, respectively.

Define $f\colon\mathcal{V}^+ \rightarrow \mathcal{V}^-$ by
$f([(z,w)],s)=([(g_\xi(z),s/w)],s)=(f_s([(z,w)]),s)$
for $([(z,w)],s)\allowbreak\in \mathcal{V}^+$. Then through the map $f$, $\mathcal{M}^+$ and $\mathcal{M}^-$ can be glued to a complex manifold~$\mathcal{M}$. Let~${\pi_{\mathcal{M}}\colon\mathcal{M} \rightarrow \hat{\Delta}}$ be the natural projection. Then $\bigl(\mathcal{M}, \hat{\Delta}, \pi_{\mathcal{M}}\bigr)$ is a complex analytic family of symmetric projective K3 surfaces over the one dimensional complex manifold~\smash{$\hat{\Delta}$} (also see \cite{Takayuki-Takato2} for a more general result).

For each fiber of the smooth compactifiable deformation $(\mathcal{F},\mathcal{S},T,\widetilde{\pi},o,j)$ constructed in Theorem~\ref{Theorem1.1}, a complex analytic family of symmetric projective K3 surfaces over one dimensional complex manifold $\hat{\Delta}=\{ s\mid s \in \Delta \backslash \{ 0 \},\, |s|<\varepsilon_0 \}$ can be constructed.

In fact, there is an analogous deformation family of quasi-projective varieties constructed here. Let $\mathcal{F}^-$ be the result manifold after blowing up $\operatorname{CP}^2\times T$ at the points
\begin{gather*}
\bigcup_{(\widetilde{\tau},\vec{x})=(\widetilde{\tau},x_1,\dots ,x_8)\in T}\{(\ell_{CP^2}\circ f_{\widetilde{\tau}}\circ f_{\widetilde{\tau},1}(x_1),\widetilde{\tau},\vec{x}),\dots ,\\
\hphantom{\bigcup_{(\widetilde{\tau},\vec{x})=(\widetilde{\tau},x_1,\dots ,x_8)\in T}} (\ell_{CP^2}\circ f_{\widetilde{\tau}}\circ f_{\widetilde{\tau},8}(x_8),\widetilde{\tau},\vec{x}),(\ell_{CP^2}\circ f_{\widetilde{\tau}}\circ f_9(\widetilde{\tau},\vec{x}),\widetilde{\tau},\vec{x}\}.
\end{gather*}
Let $\mathcal{S}^-$ be the analogous complex manifold of $\mathcal{S}$ and $\widetilde{\pi}^-\colon\mathcal{F}^- \rightarrow T^-=T$ be the map induced by $\widetilde{\pi}$. Then $(\mathcal{F}^-, T^-, \widetilde{\pi}^-)$ is a deformation family of \smash{$\operatorname{CP}^2 \#9\overline{\operatorname{CP}}^2$} analogous to the complex analytic family $(\mathcal{F}, T, \widetilde{\pi})$. Moreover, $\widetilde{\pi}^-|_{\mathcal{F}^-\backslash \mathcal{S}^-}\colon{\mathcal{F}^-\backslash \mathcal{S}^-} \rightarrow T^-$ gives a deformation family of quasi-projective varieties analogous to $(\mathcal{F},\mathcal{S},T,\widetilde{\pi},o,j)$.

In addition, let $\bigl(\mathcal{F}^+, T^+, \widetilde{\pi}^+\bigr)$ be a copy of the deformation family $(\mathcal{F}, T, \widetilde{\pi})$ and $\mathcal{S}^+$ be a copy of $\mathcal{S}$. Then the following theorem is true.

\begin{Theorem} \label{Theorem4.3}
There is a deformation family of symmetric projective K3 surfaces over a $10$-dimensional complex manifold constructed from $(\mathcal{F},\mathcal{S},T,\widetilde{\pi},o,j)$ and an analogous deformation family.
\end{Theorem}

\begin{proof}
Here, $(\mathcal{F},\mathcal{S},T,\widetilde{\pi},o,j)$ is a smooth compactifiable deformation of $S \backslash C$ over a 9-di\-men\-sion\-al complex manifold. $\widetilde{\pi}^-|_{\mathcal{F}^-\backslash \mathcal{S}^-}\colon {\mathcal{F}^-\backslash \mathcal{S}^-} \rightarrow T^-$ gives a deformation family of quasi-projective varieties analogous to $(\mathcal{F},\mathcal{S},T,\widetilde{\pi},o,j)$ and $\widetilde{\pi}^+|_{\mathcal{F}^+\backslash \mathcal{S}^+}\colon{\mathcal{F}^+\backslash \mathcal{S}^+} \rightarrow T^+$ gives a deformation family as a copy of $(\mathcal{F},\mathcal{S},T,\widetilde{\pi},o,j)$.

From the discussion above, choosing two corresponding fibers from $\mathcal{F}^+$ and $\mathcal{F}^-$, respectively, a~deformation family of symmetric projective K3 surfaces over one dimensional complex manifold~$\hat{\Delta}=\{ s\mid s \in \Delta \backslash \{ 0 \}, \, |s|<\varepsilon_0 \}$ can be constructed. Therefore, a deformation family of symmetric projective K3 surfaces over a 10-dimensional complex manifold can be constructed through the deformation families induced by $\widetilde{\pi}^-|_{\mathcal{F}^-\backslash \mathcal{S}^-}\colon{\mathcal{F}^-\backslash \mathcal{S}^-} \rightarrow T^-$ and~$\widetilde{\pi}^+|_{\mathcal{F}^+\backslash \mathcal{S}^+}\colon{\mathcal{F}^+\backslash \mathcal{S}^+} \rightarrow T^+$. So Theorem~\ref{Theorem4.3} is proved.
\end{proof}

Therefore, combining the results from Lemma~\ref{Lemma4.2} and Theorem~\ref{Theorem4.3}, the proof of Corollary~\ref{Corollary1.2} is completed.

\appendix

\section{Appendix}\label{appendixA}

The main goal of this appendix is to give the complete proofs of Lemmas~\ref{Lemma3.2} and~\ref{Lemma3.3}, Theorem~\ref{Theorem3.19} and Corollary~\ref{Corollary3.20}.

Firstly, the complete proof of Lemma~\ref{Lemma3.2} is as follows.

\begin{proof}[Proof of Lemma~\ref{Lemma3.2}]
It is obvious that \smash{$\pi_{\operatorname{CP}^2 \times U}$} is a holomorphic map. Moreover, the rank of the Jacobian matrix of \smash{$\pi_{\operatorname{CP}^2 \times U}$} is equal to 1 at every point of \smash{$\operatorname{CP}^2\times U$}. In addition, $\forall t \in U$, \smash{$\pi_{\operatorname{CP}^2 \times U}^{-1}(t)=\operatorname{CP}^2\times \{ t \}$} is a compact complex submanifold of $\operatorname{CP}^2\times U$. So $\bigl(\operatorname{CP}^2\times U, U,\pi_{\operatorname{CP}^2\times U}\bigr)$ is a complex analytic family of $\operatorname{CP}^2$ over a 1-dimensional complex manifold $U$.

Now \smash{$\widetilde{\mathcal{S}}=\bigl\{ (x, \tau) \in \operatorname{CP}^2 \times U|x \in f_{\tau}(C_0(\tau)) \bigr\}$} is a complex submanifold of $\operatorname{CP}^2 \times U$ and the rank of the Jacobian matrix of \smash{$\pi_{\operatorname{CP}^2 \times U}|_{\widetilde{\mathcal{S}}}$} is equal to 1 at every point of \smash{$\widetilde{\mathcal{S}}$}. $\forall t \in U$, $\smash{\pi_{\operatorname{CP}^2 \times U}|_{\widetilde{\mathcal{S}}}^{-1}(t)}=f_{t}(C_0(t))\times \{ t \}$ is a compact complex submanifold of $\widetilde{\mathcal{S}}$. So \smash{$\bigl(\widetilde{\mathcal{S}},U,\pi_{\operatorname{CP}^2 \times U}|_{\widetilde{\mathcal{S}}}\bigr)$} is a complex analytic family of elliptic curves over the 1-dimensional complex manifold $U$.

Therefore, $\bigl(\widetilde{\mathcal{S}},U,\pi_{\operatorname{CP}^2 \times U}|_{\widetilde{\mathcal{S}}}\bigr)$ is a complex analytic family of elliptic curves and $\bigl(\operatorname{CP}^2\times U, U,$ $\pi_{\operatorname{CP}^2\times U}\bigr)$ is a complex analytic family of $\operatorname{CP}^2$. So Lemma~\ref{Lemma3.2} is proved. \end{proof}

Secondly, the proof of Lemma~\ref{Lemma3.3} was simply described in \cite{Takayuki-Takato1}. Here, a complete proof of Lemma~\ref{Lemma3.3} is as follows.

\begin{proof}[Proof of Lemma~\ref{Lemma3.3}]
For all $\widetilde{\tau} \in U$, let $N_{\widetilde{\tau}}$ be the normal bundle of the strict transform for the elliptic curve $f_{\widetilde{\tau}}(C_0(\widetilde{\tau}))\times \{ \widetilde{\tau} \} \cong C_0(\widetilde{\tau})$ after blowing up at nine points
\[\{ \widetilde{p}_1=(f_{\widetilde{\tau}}(p_1),\widetilde{\tau}),\dots , \widetilde{p}_9=(f_{\widetilde{\tau}}(p_9),\widetilde{\tau}) \}\subset f_{\widetilde{\tau}}(C_0(\widetilde{\tau}))\times \{ \widetilde{\tau} \}\]
with
\[ \{ p_1=f_{\widetilde{\tau},1}(\hat{p}_1), p_2=f_{\widetilde{\tau},2}(\hat{p}_2),\dots ,p_8=f_{\widetilde{\tau},8}(\hat{p}_8),p_9 \}\subset C_0(\widetilde{\tau})=\mathbb{C}/ \langle 1,\widetilde{\tau} \rangle\]
and $(\hat{p}_1,\dots ,\hat{p}_8)\in $$U_1 \times \dots \times U_8$. In addition, let $C(\widetilde{\tau})\times \{ \widetilde{\tau} \}$ be the strict transform of $ f_{\widetilde{\tau}}(C_0(\widetilde{\tau}))\times \{ \widetilde{\tau} \}$. Here, $N_{\widetilde{\tau}}$ is assumed to satisfy the Diophantine condition.

Then $N_{\widetilde{\tau}}$ is isomorphic to
\[
\mathcal{O}_{\operatorname{CP}^2}(3)|_{f_{\widetilde{\tau}}(C_0(\widetilde{\tau}))}\otimes \mathcal{O}_{f_{\widetilde{\tau}}(C_0(\widetilde{\tau}))}(-f_{\widetilde{\tau}}(p_1)-f_{\widetilde{\tau}}(p_2)-\dots -f_{\widetilde{\tau}}(p_9))\in \operatorname{Pic}^0(f_{\widetilde{\tau}}(C_0(\widetilde{\tau}))).
\]
Furthermore, taking use of notations from Section~\ref{sec2}, since $V_s$ is biholomorphic to a topologically trivial annulus bundle over the strict transform $C(\widetilde{\tau})$, $V_s$ is homotopic to $ S_1 \times S_2 \times S_3$ with~$ { \{S_1 , S_2 , S_3 \}}$ being circles and $ S_1 \times S_2$ being $C^\infty$ sections for $s\in \Delta \setminus \{ 0\}$ and sufficiently small \cite{Takayuki-Takato2} (also see Theorem~\ref{Theorem2.8}). Then through integrating the global holomorphic form $\sigma_s$ defined in Section~\ref{sec2} on $ S_1 \times S_2$, the result value is $q-p \cdot \widetilde{\tau}$ \cite{Takayuki-Takato1}. So there is an equation
\[9p_0-\sum _{j=1}^8p_j-p_9=q-p\cdot \widetilde{\tau} \quad \mod \langle 1, \widetilde{\tau} \rangle \]
with $p_0$ being an inflection point on $C_0(\widetilde{\tau})$ derived in \cite{Takayuki-Takato2}. Then for $\hat{t}=(\widetilde{\tau},\hat{p}_1,\dots ,\hat{p}_8)\in T$, $p_9=f_9\bigl(\hat{t}\bigr)\in C_0(\widetilde{\tau})=\mathbb{C}/ \langle 1, \widetilde{\tau} \rangle$.

Therefore, $N_{\widetilde{\tau}}$ satisfies the Diophantine condition corresponding to $(p,q)$ if and only if
\[9p_0-\sum _{j=1}^8p_j-p_9=q-p\cdot \widetilde{\tau} \quad \mod \langle 1, \widetilde{\tau} \rangle \]
with $p_0$ being an inflection point on $C_0(\widetilde{\tau})$.

In addition, taking use of the detailed definition of the map $\widetilde{\pi}$, $(\mathcal{F},T, \widetilde{\pi})$ and $(\mathcal{S},T, \widetilde{\pi}|_\mathcal{S})$ can be proved to be complex analytic families directly.

In conclusion, through blowing up nine points on every fiber of $\bigl(\operatorname{CP}^2 \times T,T, \pi\bigr)$ such that the normal bundle of the strict transform of the elliptic curve mentioned in Section~\ref{sec3.1} embedded in each fiber satisfies the Diophantine condition corresponding to the Diophantine number pair~$(p,q)$, the complex analytic families $(\mathcal{F},T, \widetilde{\pi})$ and $(\mathcal{S},T, \widetilde{\pi}|_\mathcal{S})$ can be constructed. So Lemma~\ref{Lemma3.3} is proved.
\end{proof}

Thirdly, the complete proof of Theorem~\ref{Theorem3.19} is as follows.

\begin{proof}[Proof of Theorem~\ref{Theorem3.19}]
First of all, it is to prove that there exists $b>0$ such that the metric $h \cdot {\rm e}^{-b \Psi_s}$ is a~smooth hermitian metric on $L|_{S(\tau)\backslash C(\tau)}$ with $\Theta_{h \cdot {\rm e}^{-b \Psi_s}}>0$.

On the set $S(\tau)\backslash W$, since $\Psi_s$ is constant, $\partial \bar{\partial}\Psi_s=0$. Then
\[\Theta_{h \cdot {\rm e}^{-b \Psi_s}}=\Theta_h + \sqrt{-1} \frac{b}{\pi} \partial \bar{\partial} \Psi_s =\Theta_h=\Theta_{h_L}>0.\]

On the set $\{ [(z,w)]\in W\backslash C(\tau) \mid r_2 \leq |w|<r \}$, from the construction of $\Psi_s$ above, since~$\widetilde{\theta}_s(x)$ is smooth for $x>0$, the coefficients of $\partial \bar{\partial}\Psi_s$ are bounded. Moreover, since $\Theta_{h_L}>0$ on~$S(\tau)\backslash C(\tau)$, the coefficients of $\Theta_h$ have a positive lower bound on $\{ [(z,w)]\in W\backslash C(\tau) \mid r_2 \leq |w|<r \}$. Then $\exists$ sufficiently small $b>0$ such that $\Theta_{h \cdot {\rm e}^{-b \Psi_s}}=\Theta_h + \sqrt{-1} \frac{b}{\pi}\partial \bar{\partial} \Psi_s >0$.

On the set $\{ [(z,w)]\in W\backslash C(\tau) \mid |w|<r_2 \}$,
$\sqrt{-1} \partial \bar{\partial} \Psi_s=2 \sqrt{-1} \cdot {\rm d}w \wedge {\rm d}\bar{w} /|w|^2>0$.
So $\Theta_{h \cdot {\rm e}^{-b \Psi_s}}=\Theta_h + \sqrt{-1} \frac{b}{\pi}\partial \bar{\partial} \Psi_s >0$.

The next step is to prove that $\Theta_{h \cdot {\rm e}^{-b \Psi_s}}$ is of the form
\[
\omega|_{W_{\varepsilon_0}}=\frac{ (L.C(\tau))}{\operatorname{ Im}\tau}\cdot \sqrt{-1} {\rm d}z \wedge {\rm d} \bar{z}+\frac{2b}{\pi} \cdot \frac{\sqrt{-1}{\rm d}w \wedge {\rm d} \bar{w}}{|w|^2}
\]
on the set $W_{\varepsilon_0}=\{[(z,w)]\in W \backslash C(\tau)\mid |w|<\sqrt{\varepsilon_0}r\}$.

On the set $W_{\varepsilon_0}=\{[(z,w)]\in W \backslash C(\tau)\mid |w|<\sqrt{\varepsilon_0}r\}$,
\begin{align*}
\Theta_{h \cdot {\rm e}^{-b \Psi_s}}&=\Theta_{(\pi_W)^*h_{C(\tau)}}+\frac{2b}{\pi}\cdot \frac{\sqrt{-1}{\rm d}w \wedge {\rm d} \bar{w}}{|w|^2}\\
&=\frac{(L.C(\tau))}{\operatorname{ Im}\tau}\cdot \sqrt{-1} {\rm d}z \wedge {\rm d} \bar{z}+\frac{2b}{\pi}\cdot \frac{\sqrt{-1}{\rm d}w \wedge {\rm d} \bar{w}}{|w|^2}.
\end{align*}
 So the curvature form $\Theta_{h \cdot {\rm e}^{-b\Psi_s}}$ gives a complete K\"ahler metric on $S(\tau)\backslash C(\tau)$.

The last step is to prove that $\omega|_{W_{\varepsilon_0}}$ is Ricci-flat on the set $W_{\varepsilon_0}=\{ [(z,w)]\in W \backslash C(\tau)\mid |w|<\sqrt{\varepsilon_0}r \}$.

The coefficient matrix of $\omega|_{W_{\varepsilon_0}}$ on the set $W_{\varepsilon_0}=\{ [(z,w)]\in W \backslash C(\tau)\mid |w|<\sqrt{\varepsilon_0}r \}$ is \[(g_{ij})_{2 \times 2}= \begin{pmatrix}
 \frac{2 (L.C(\tau))}{\operatorname{ Im}\tau} &0\\
 0 & \frac{4b}{\pi |w|^2}
\end{pmatrix}.\]
Then the determinate of the coefficient matrix is
\[\det((g_{ij})_{2\times 2})=\frac{8b (L.C(\tau))}{\pi \operatorname{ Im}\tau \cdot |w|^2}.\]
So $\operatorname{Ric}(\omega|_{W_{\varepsilon_0}})=0$. Therefore, $\omega|_{W_{\varepsilon_0}}$ is Ricci-flat.

In conclusion, $\exists \ b>0$ such that $\Theta_{h\cdot {\rm e}^{-b\Psi_s}}$ gives a complete K\"ahler metric on $S(\tau)\backslash C(\tau)$. In addition, $\Theta_{h\cdot {\rm e}^{-b\Psi_s}}$ is of the form
\[
\omega|_{W_{\varepsilon_0}}=\frac{(L.C(\tau))}{\operatorname{ Im}\tau}\cdot \sqrt{-1} {\rm d}z \wedge {\rm d} \bar{z}+\frac{2b}{\pi} \cdot \frac{\sqrt{-1}{\rm d}w \wedge {\rm d} \bar{w}}{|w|^2}
\]
on the set $W_{\varepsilon_0}=\{[(z,w)]\in W \backslash C(\tau)\mid |w|<\sqrt{\varepsilon_0}r\}$. Furthermore, $\omega|_{W_{\varepsilon_0}}$is Ricci-flat. That is to say, $\omega|_{W_{\varepsilon_0}}$ corresponds to a Ricci-flat K\"ahler metric. So Theorem~\ref{Theorem3.19} is proved.
\end{proof}

Finally, the complete proof of Corollary~\ref{Corollary3.20} is as follows.

\begin{proof}[Proof of Corollary~\ref{Corollary3.20}]
Since $\widetilde{\pi}|_{\mathcal{F}\backslash \mathcal{S}}$ is a holomorphic surjective submersion and $f$ is a biholomorphic map from~$D_3$ to $f(D_3)$, then $f \circ \widetilde{\pi}|_{\mathcal{F}\backslash \mathcal{S}}$ is also a holomorphic surjective submersion. $(f \circ \widetilde{\pi}|_{\mathcal{F}\backslash \mathcal{S}})^{-1}(0)=j(S\backslash C)$.

Moreover, since $\widetilde{\pi}|_{\mathcal{F}\backslash \mathcal{S}}$ is locally trivial in the $C^\infty$ category and $f$ is a diffeomorphism, then \smash{$f \circ \widetilde{\pi}|_{\mathcal{F}\backslash \mathcal{S}}$} is also locally trivial in the $C^\infty$ category.

In conclusion, $f \circ \widetilde{\pi}|_{\mathcal{F}\backslash \mathcal{S}}\colon \mathcal{F} \rightarrow f(D_3)$ is a deformation family of the non-compact complete K\"ahler manifold $j(S\backslash C)$ over a 9-dimensional complex manifold $f(D_3)$ satisfying the conditions in Definition~\ref{Definition3.5}. So Corollary~\ref{Corollary3.20} is proved.
\end{proof}

\subsection*{Acknowledgements}
Thanks Professor Zhou Zhang for his course on algebraic geometry. Thanks Professor Sho Tanimoto for providing some useful materials. In particular, thanks Professor Ryoichi Kobayashi for introducing this topic and giving invaluable comments. I would also like to thank the anonymous referees for their suggestions and comments.

\pdfbookmark[1]{References}{ref}
\LastPageEnding


\begin{thebibliography}{99}
\footnotesize\itemsep=0pt

\bibitem{Anna}
Abasheva A., D\'eev R., Complex surfaces with many algebraic structures,
 \href{https://doi.org/10.1093/imrn/rnad190}{\textit{Int. Math. Res. Not.}}
 \textbf{2024} (2024), 7379--7400,
 \href{http://arxiv.org/abs/2303.10764}{arXiv:2303.10764}.

\bibitem{Yukitaka-Klaus}
Abe Y., Kopfermann K., Toroidal groups. Line bundles, cohomology and
 quasi-abelian varieties, \textit{Lecture Notes in Math.}, Vol.~1759,
 \href{https://doi.org/10.1007/b80605}{Springer}, Berlin, 2001.

\bibitem{Vladimir}
Arnold V.I., Bifurcations of invariant manifolds of differential equations, and
 normal forms of neighborhoods of elliptic curves,
 \href{https://doi.org/10.1007/BF01076024}{\textit{Funct. Anal. Appl.}}
 \textbf{10} (1976), 249--259.

\bibitem{Edoardo-Elizabeth-Francisco}
Ballico E., Gasparim E., Rubilar F., 20~open questions about deformations of
 compactifiable manifolds,
 \href{https://doi.org/10.1007/s40863-021-00213-8}{\textit{S\~{a}o Paulo~J.
 Math. Sci.}} \textbf{15} (2021), 661--681,
 \href{http://arxiv.org/abs/2004.11299}{arXiv:2004.11299}.

\bibitem{Marco}
Brunella M., On {K}\"ahler surfaces with semipositive {R}icci curvature,
 \textit{Riv. Math. Univ. Parma~(N.S.)} \textbf{1} (2010), 441--450.

\bibitem{Jean}
Demailly J.-P., Complex analytic and differential geometry, Universit\`e de
 Grenoble~I, Grenoble, 1997.

\bibitem{Lawrence}
Evans L.C., Partial differential equations, \textit{Grad. Stud. Math.},
 Vol.~19, \href{https://doi.org/10.1090/gsm/019}{American Mathematical
 Society}, Providence, RI, 1998.

\bibitem{Yoshio}
Fujimoto Y., On rational elliptic surfaces with multiple fibers,
 \href{https://doi.org/10.2977/prims/1195171661}{\textit{Publ. Res. Inst.
 Math. Sci.}} \textbf{26} (1990), 1--13.

\bibitem{Elizabeth-Francisco}
Gasparim E., Rubilar F., Deformations of noncompact {C}alabi-{Y}au manifolds,
 families and diamonds, in Geometry at the {F}rontier---{S}ymmetries and
 {M}oduli {S}paces of {A}lgebraic {V}arieties, \textit{Contemp. Math.}, Vol.~766, \href{https://doi.org/10.1090/conm/766/15378}{American Mathematical
 Society}, Providence, RI, 2021, 117--132,
 \href{http://arxiv.org/abs/1908.09045}{arXiv:1908.09045}.

\bibitem{Phillip}
Griffiths P., Harris J., Principles of algebraic geometry, \textit{Wiley Classics Lib.},
 \href{https://doi.org/10.1002/9781118032527}{John Wiley \& Sons}, New York,
 1994.

\bibitem{Vincent-Ahmed}
Guedj V., Zeriahi A., Degenerate complex {M}onge--{A}mp\`ere equations,
 \textit{EMS Tracts Math.}, Vol.~26,
 \href{https://doi.org/10.4171/167}{European Mathematical Society (EMS)},
 Z\"urich, 2017.

\bibitem{Robin}
Hartshorne R., Algebraic geometry, \textit{Grad. Texts in Math.}, Vol.~52,
 \href{https://doi.org/10.1007/978-1-4757-3849-0}{Springer}, New York, 1977.

\bibitem{Hans}
Hein H.J., Complete Calabi--Yau metrics
 from~$\mathbb{P}^2\#9\overline{\mathbb{P}}^2$,
 \href{http://arxiv.org/abs/1003.2646}{arXiv:1003.2646}.

\bibitem{Kunihiko}
Kodaira K., Complex manifolds and deformation of complex structures, \textit{Class.
 Math.}, \href{https://doi.org/10.1007/b138372}{Springer}, Berlin, 2005.

\bibitem{Takayuki-Takato1}
Koike T., Uehara T., A gluing construction of projective~{K}3 surfaces,
 \href{http://arxiv.org/abs/1903.01444}{arXiv:1903.01444}.

\bibitem{Takayuki-Takato2}
Koike T., Uehara T., A gluing construction of {K}3 surfaces,
 \href{https://doi.org/10.46298/epiga.2022.volume6.8504}{\textit{\'Epijournal
 G\'eom. Alg\'ebrique}} \textbf{6} (2022), 12, 15~pages,
 \href{http://arxiv.org/abs/2108.07168}{arXiv:2108.07168}.

\bibitem{Julius-David}
Ross J., Nystr\"om D.W., Homogeneous {M}onge--{A}mp\`ere equations and
 canonical tubular neighbourhoods in {K}\"ahler geometry,
 \href{https://doi.org/10.1093/imrn/rnw200}{\textit{Int. Math. Res. Not.}}
 \textbf{2017} (2017), 7069--7108,
 \href{http://arxiv.org/abs/1403.3282}{arXiv:1403.3282}.

\bibitem{Alexandre}
Sukhov A., Regularized maximum of strictly plurisubharmonic functions on an
 almost complex manifold,
 \href{https://doi.org/10.1142/S0129167X13500973}{\textit{Internat.~J. Math.}}
 \textbf{24} (2013), 1350097, 6~pages,
 \href{http://arxiv.org/abs/1303.5312}{arXiv:1303.5312}.

\bibitem{Tomasz-Halszka}
Szemberg T., Tutaj-Gasi\'nska H., General blow-ups of the projective plane,
 \href{https://doi.org/10.1090/S0002-9939-02-06488-2}{\textit{Proc. Amer.
 Math. Soc.}} \textbf{130} (2002), 2515--2524.

\bibitem{Tetsuo}
Ueda T., On the neighborhood of a compact complex curve with topologically
 trivial normal bundle,
 \href{https://doi.org/10.1215/kjm/1250521670}{\textit{J.~Math. Kyoto Univ.}}
 \textbf{22} (1982), 583--607.

\bibitem{Christian}
Vogt C., Line bundles on toroidal groups,
 \href{https://doi.org/10.1515/crll.1982.335.197}{\textit{J.~Reine Angew.
 Math.}} \textbf{335} (1982), 197--215.

\bibitem{Claire}
Voisin C., Hodge theory and complex algebraic geometry.~{I}, \textit{Cambridge
 Stud. Adv. Math.}, Vol.~76,
 \href{https://doi.org/10.1017/CBO9780511615344}{Cambridge University Pres}s,
 Cambridge, 2002.

\end{thebibliography}
\end{document}